\documentclass[11pt,a4paper]{article} 

\usepackage[T1]{fontenc}
\usepackage[utf8]{inputenc}

\usepackage{cite}

\usepackage{amssymb}
\usepackage{latexsym}
\usepackage{amsmath}
\usepackage{mathtools}
\usepackage{amscd}
\usepackage{graphicx}

\usepackage{proof}

\usepackage[hmargin=3cm,vmargin=3.5cm]{geometry}

\newtheorem{theorem}{Theorem}[section]

\newtheorem{definition}[theorem]{Definition}

\newtheorem{example}[theorem]{Example}

\newtheorem{remark}[theorem]{Remark}

\newcommand{\dontshow}[1]{}

\usepackage{MnSymbol,wasysym}
\usepackage{newunicodechar}

\newunicodechar{〈}{\ensuremath{\langle}}
\newunicodechar{〉}{\ensuremath{\rangle}}
\newunicodechar{×}{\ensuremath{\times}}
\newunicodechar{∈}{\ensuremath{\in}}
 \newunicodechar{≐}{\ensuremath{\dot{=}}}
\newunicodechar{∧}{\ensuremath{\land}}
\newunicodechar{∨}{\ensuremath{\lor}}
\newunicodechar{∃}{\ensuremath{\exists}}
\newunicodechar{∀}{\ensuremath{\forall}}
\newunicodechar{λ}{\ensuremath{\lambda}}
\newunicodechar{≈}{\ensuremath{\approx}}
\newunicodechar{→}{\ensuremath{\rightarrow}}
\newunicodechar{↔}{\ensuremath{\leftrightarrow}}

\newunicodechar{Σ}{\ensuremath{\Sigma}}
\newunicodechar{Π}{\ensuremath{\Pi}}
\newunicodechar{σ}{\ensuremath{\sigma}}
\newunicodechar{π}{\ensuremath{\pi}}

\newunicodechar{Γ}{\ensuremath{\Gamma}}
\newunicodechar{Δ}{\ensuremath{\Delta}}
\newunicodechar{Θ}{\ensuremath{\Theta}}

\newunicodechar{★}{\ensuremath{\star}}
\newunicodechar{•}{\ensuremath{\cdot}}
\newunicodechar{⁻}{\ensuremath{{}^{-}}}
\newunicodechar{¹}{\ensuremath{{}^{1}}}
\newunicodechar{⇒}{\ensuremath{\Rightarrow}}
\newunicodechar{ °}{\ensuremath{\circ}}
\newunicodechar{ α}{\ensuremath{\alpha}}
\newunicodechar{ β}{\ensuremath{\beta}}
\newunicodechar{‣}{\ensuremath{\blacktriangleright}}
\newunicodechar{⇛}{\ensuremath{\Rrightarrow}}
\newunicodechar{⋃}{\ensuremath{\union}}
\newunicodechar{⊙}{\ensuremath{\odot}}
\newunicodechar{⊗}{\ensuremath{\otimes}}
\newunicodechar{⊕}{\ensuremath{\oplus}}
 
\newunicodechar{κ}{\ensuremath{\kappa}}
\newunicodechar{°}{\ensuremath{{}^{\circ}}}
  
\newunicodechar{⌢}{\ensuremath{\frown}}
\newunicodechar{↓}{\ensuremath{\downarrow}}
\newunicodechar{↑}{\ensuremath{\uparrow}}
\newunicodechar{▷}{\ensuremath{\rhd}}
\newunicodechar{₀}{\ensuremath{{}_0}}
\newunicodechar{₁}{\ensuremath{{}_1}}
\newunicodechar{∼}{\ensuremath{\sim}}
\newunicodechar{≤}{\ensuremath{\le}}

\input diagxy.tex
\xyoption{curve}

\begin{document}

\newcommand{\sbrac}[1]{[\![#1]\!]}
\newcommand{\Hom}{{\rm Hom}}
\newcommand{\True}{\top}
\newcommand{\False}{\bot}
\newcommand{\drule}[2]{\frac{#1}{#2}}
\newcommand{\pgate}{\vdash}
\newcommand{\sequent}{\Longrightarrow}
\newcommand{\sequentunder}[1]{\xRightarrow{\;\;#1\;\;}}

\newcommand{\nat}{{\mathbb N}}
\newcommand{\nattype}[0]{{\rm N}}
\newcommand{\context}{{\rm context}}
\newcommand{\type}{{\rm type}}
\newcommand{\form}{{\rm formula}}

\newcommand{\Ty}{{\rm Ty}}
\newcommand{\Tm}{{\rm Tm}}

\newcommand{\app}{{\rm app}}

\newcommand{\vf}[0]{{\varphi}}
\newcommand{\pto}[0]{\rightharpoondown}
\newcommand{\piev}[0]{{\rm ev}}
\newcommand{\memof}[0]{\, \epsilon \,}
\newcommand{\subseteqof}[0]{\, \dot{\subseteq} \,}

\newcommand{\mono}[0]{\to/ >->/}

\newcommand{\emptycontext}{\langle\rangle}
\newcommand{\emptyctx}{\langle\rangle}
\newcommand{\codenote}[1]{\footnote{\tt #1}}

\title{From type theory to setoids and back} 
\author{Erik Palmgren \\ Department of Mathematics, Stockholm University}
\date{September 16, 2019}

\maketitle

\begin{abstract} 
A model of Martin-Löf extensional type theory with universes is formalized in Agda, an interactive proof system based on Martin-Löf intensional type theory.
This may be understood, we claim, as a solution to the old problem of modeling the full extensional theory
in the intensional theory. Types are interpreted as setoids, and the model is therefore a setoid model. We solve the problem of intepreting type universes by utilizing Aczel's type of iterative sets, and show how it can be made into a setoid of small setoids containing the necessary setoid constructions.

In addition we interpret the bracket types of Awodey and Bauer. Further quotient types should be interpretable.

{\flushleft AMS} Mathematical Subjects Classification : 03B15, 03B35, 03E70, 03F50 

\end{abstract}

\tableofcontents

\section{Introduction}

In this paper we present an interpretation of full extensional Martin-Löf type theory \cite{martinlof1982} into intensional Martin-Löf type theory via setoid constructions. There are several qualifications to this statement. We actually formalise this interpretation in the Agda proof assistant, where the fragment of the system used is  considered as the intensional type theory. Our system of extensional type theory is not a syntactically defined system, but rather a system of closure rules for judgements about setoids in Agda. The fragment of Agda used is limited to certain kinds of inductive-recursive definitions and record types. The K rule of Agda is not used.  We believe that the proofs carried out in this fragment may also be carried out in a Logical Framework presentation of Martin-Löf type theory with a {\em super universe}\cite{palmgren1998}. A super universe is closed under construction of universes, in addition to the standard type constructions introduced in \cite{martinlof1982}.

A first approach that may come to mind when interpreting extensional type theory using setoid constructions
 is to interpret a family of types over a context as a family of setoids over a setoid that interprets the context.  (For background on setoids see Section \ref{prelsec}.)
 The basic judgement forms of Martin-L\"of type theory \cite{martinlof1984} are displayed to the left in the table below. 
$$
\begin{array}{llll}
&\Gamma \Longrightarrow A \; {\rm type}  &\qquad & A  :{\rm Fam}(\Gamma) \\
&\Gamma \Longrightarrow A  = B  && ? \\
&\Gamma \Longrightarrow a  : A && a : \Pi(\Gamma, A) \\
& \Gamma \Longrightarrow a = b : A && a =_{\Pi(\Gamma, A)} b 
\end{array}
$$
We may now try to interpret the forms of judgements as the statements about setoids to the right. But we do not yet have any obvious interpretation of the type equality.
We need to compare $A$ and $B$ as setoid families over $\Gamma$. A crucial problem is how to interpret the type equality rule 
\begin{equation} \label{tyeq}
\infer{\Gamma \Longrightarrow a  :  B}{\Gamma \Longrightarrow a  :  A & \Gamma \Longrightarrow A  = B}
\end{equation}

A solution is to embed all dependent families of setoids in to a big universal setoid (or as we will call it {\em classoid}). To obtain a setoid model without coherence problems we may seek inspiration from  type-free interpretations of (extensional) type theory; see Aczel \cite{aczel1977}, Smith \cite{smith1984}, Beeson \cite[Ch.\ XI]{beeson1985}. But instead of using combinators or recursive realizers as type free objects, we use constructive iterative sets in the sense of Aczel \cite{aczel1978}.

Aczel's type of iterative sets $V$ \cite{aczel1978} consists of well-founded trees where the branching $f$ can be indexed by any type $A$ in a universe $U$ of small types. The introduction rule tells how to build a set $\alpha={\rm sup}(A,f)$ from a family $f(x)$ ($x: A$) of previously constructed sets

$$\infer[(V\; {\rm intro})]{{\rm sup}(A,f) : V}{ A : U & f: A \to V}$$
Equality $=_V$ is defined by the smallest bisimulation, and then membership is given by
$$x \in_V {\rm sup}(A,f) := (\exists a: A) (x =_V f(a)).$$
The classoid $ {\mathbb V} = (V,=_V)$ forms, together with the membership relation $\in_V$, a model of Constructive Zermelo-Fraenkel set theory (CZF) with Dependent Choice (DC), and possibly further axioms, depending on the type theory. It is thus expected to be a rich universe of sets. In fact, each set $\alpha= {\rm sup}(A,f)$ may also be understood as a setoid on the type $A$
\begin{equation} \label{can-std}
\kappa(\alpha) = (A, =_f)
\end{equation}
where $a =_f b$ is defined as $ f(a) =_V f(b)$. 
The assignment $\kappa$ may be extended to a full and faithful functor from the category of sets in $V$ to the category of small setoids. Using $\kappa$ we can also construct a bijection of classoids
$$  {\mathbb V}  \cong {\rm Sub}({\mathbb V}).$$
Following Aczel \cite{aczel1982} one can see that there are internal versions in $V$ of the setoid construction for $\Pi$, $\Sigma$ and extensional identity, which commute with $\kappa$
$$ \kappa(\sigma(\alpha,f)) \cong \Sigma(\kappa(\alpha), \kappa \circ f) \qquad
\kappa(\pi(\alpha,f)) \cong \Pi(\kappa(\alpha), \kappa \circ f)$$
Thus ${\mathbb V}$ is suitable for interpreting both terms and types of dependent type theory.
A type $A$ in a context $\Gamma$ will be interpreted as an extensional function $A : \kappa(\Gamma) \rightarrow {\mathbb V}$. Now any two types in the same context can be compared. A raw term $a$ in a context $\Gamma$ will likewise be interpreted as extensional function $a : \kappa(\Gamma) \rightarrow {\mathbb V}$. 
The judgement $a : A$ will be interpreted as membership judgement. The new {\em setoid interpretation} is on the right in the table below.
\begin{equation} \label{newinterp}
\begin{array}{llll}
&\Gamma \Longrightarrow A \; {\rm type}  &\qquad &   A : [\kappa(\Gamma) \rightarrow {\mathbb V}]\\
&\Gamma \Longrightarrow A  = B  &&   A =_{\rm ext} B: [\kappa(\Gamma) \rightarrow  {\mathbb V}]\\
&\Gamma \Longrightarrow a  :  A &&  \forall x : \kappa(\Gamma), a(x) \in_V A(x) \\
& \Gamma \Longrightarrow a = b : A &&   \forall x : \kappa(\Gamma), a(x)=_V b(x) \in_V A(x)
\end{array}
\end{equation}
The interpretation of the problematic type equality rule (\ref{tyeq}) is now direct. Further the basic rules in type theory for $\Sigma$, $\Pi$, $+$, extensional identity types, and the basic types $N_0$ and $N$  can now be interpreted; see Section \ref{ml0}. Some further considerations are necessary to interpret the hierarchy of type universes. Here we use a superuniverse \cite{palmgren1998, rathjen2000, rathjen2001}, which is a type universe closed under the operation of building a universe over a family of base types. This makes it possible to build the hierarchy of setoid universes internally to the superuniverse, and interpret the universe rules à la Russell. This is covered in Section \ref{mlomega}. Bracket type constructions are defined in Section \ref{brack}. The interpretation of the judgement forms is fixed in Section \ref{interpretation}. Section \ref{listrules} lists all the rules interpreted together with references to the formalization. Section \ref{overview} contains links to the actual formalization which is available on-line. A comparison between Agda and the Logical Framework is made in Section \ref{LFAgda}.

\subsection*{Related work} 
Though it was widely recognized from the beginning of Martin-Löf type theory that it had a natural classical set-theoretic interpretation, a proof seems only to have been written down and published in detail by Salvesen in her 1986 MSc thesis\cite{salvesen}.  Closely related to the work of the present paper is that of Aczel \cite{aczel1999} who interprets extensional Martin-Löf type theory with universes in an extension of CZF with a hierarchy of inaccessible sets. Earlier Werner \cite{werner1997} had modelled a Coq system in ZFC and vice versa ZFC in Coq using Aczel's encoding of sets. A refinement by Barras models a Coq system in {\em intuitionistic} ZF, and formalizes the model in Coq \cite{barras2010}.  Rathjen and Tupailo \cite{rathjen+tupailo} make a close analysis of the interpretation of CZF into Martin-Löf type theory, and the question about what general classes of set-theoretic statements are validated in type theory.

Hofmann \cite{hofmann1997a} modelled an extensional Martin-Löf type theory ${\rm TT}_{\rm E}$ (without universes) in an "intensional" version of the theory ${\rm TT}_{\rm I}$.  A conservativity theorem (for type inhabitation) of ${\rm TT}_{\rm E}$ over ${\rm TT}_{\rm I}$ is established \cite[Thm 3.2.5]{hofmann1997a}.  Note that ${\rm TT}_{\rm I}$ has function extensionality and the UIP axiom, so it is different from what is usually called {\em intensional Martin-Löf type theory} \cite{nordstrom+peterson+smith}. Hofmann  \cite[Ch 5.1, 5.3]{hofmann1997a} also constructs setoid models of quotient types. The setoids are Prop-valued, or proof-irrelevant, in the sense that the truth values of equalities are in the type of propositions ${\rm Prop}$ rather than in ${\rm Set}$, the type of small types, as in the Definition \ref{setoid} we use below. The equality of morphisms and families in these setoid models are however definitional, which is an another difference.

The Minimalist Foundation of Maietti and Sambin \cite{maietti+sambin} is a two level type theory consisting of an extensional theory and a more fundamental theory, Minimal type theory, which is intensional. The extensional level is modelled \cite{maietti2009} into the minimal theory using a quotient construction.  Moreover equality of morphisms in this model is extensional as is ours. Further related constructions of extensional and quotients structures are  in \cite{maietti+rosolini2012, maietti+rosolini2013}.

Some obstacles and opportunities for constructing convenient categories and universes of setoids inside intensional Martin-Löf type theory are demonstrated and discussed in \cite{wilander2010, wilander2012, palmgren2012, palmgren+wilander, palmgen2018}.

\subsection*{Acknowledgements}
The author is grateful to the Hausdorff Research Institute for Mathematics in Bonn for the invitation to the trimester program Constructions, Sets and Types during Summer 2018, and to the scientific organizing committee:  Douglas S. Bridges, Michael Rathjen, Peter Schuster and Helmut Schwichtenberg. A first version of this work was developed and presented there. The author also wants to thank Thierry Coquand,  Giovanni Sambin, and especially, Maria Emilia Maietti, for providing references and comments. Thanks to Peter Dybjer for an invitation the workshop {\em Logic and Types 2018} in G\"oteborg to give a presentation of this work. 

 \section{Preliminaries}  \label{prelsec}

We recall some definitions and facts around setoids. Note that we use the propositions-as-types principle throughout. In particular the notion of setoid uses this principle (in contrast to e.g.\ setoids of the standard library in the Coq system).

\subsection{Setoids}

\begin{definition} \label{setoid}
 {\em
 A {\em setoid } $A=(|A|,=_A)$ is type $A$ together with an equivalence relation $=_A$, so that $x=_Ay$ is type.  A {\em setoid map, or extensional function} $f: A \rightarrow B$ is a pair $f=(|f|, {\rm ext}_f)$ consisting of a function $|f| : |A| \rightarrow |B|$ and ${\rm ext}_f$ a proof of extensionality, i.e. that
$$(\forall x,y : |A|)[x =_A y \Rightarrow |f|(x) =_B |f|(y)].$$
}
\end{definition}
Write $a : A$ for $a : |A|$, and $f(x) = |f|(x)$ for a setoid $A$ and a setoid map $f$.

\medskip
For setoids $A$ and $B$ the {\em product setoid} $A\times B$ 
is given by
$$|A\times B| =_{\rm def} |A| \times |B|$$
and
$$(a,b) =_{A \times B} (c,d) \Longleftrightarrow_{\rm def} a=_A c \land b =_B d.$$
For setoids $A$ and $B$ the {\em exponent setoid} $[A \rightarrow B] = B^A $ 
is given by
$$|B^A| =_{\rm def} (\Sigma f: |A| \rightarrow |B|)(\forall x,y: |A|)(x=_A y \Rightarrow f(x) =_B f(y))$$
and
$$(f,p) =_{B^A} (g,q) \Longleftrightarrow_{\rm def} (\forall x: |A|)(f(x)=_B g(x)).$$

With type universes we may introduce some distinctions of setoids which are useful and necessary to solve predicativity problems. Let $U_n$, $n=0,1,2,\ldots$ denote the cumulative universes of a Martin-L\"of type theory, \`a la Russell. (In Agda and Coq these are available as Set0, Set1, Set2, ...  resp Type0, Type1, Type2, ... .)  We recall a definition and some examples from \cite{palmgen2018b}:

\begin{definition} {\em
An {\em $(m, n)$-setoid $A=(|A|,=_A)$} is a type $|A| \in U_m$ with an equivalence relation $=_A : |A| \rightarrow |A| \rightarrow U_n$. An $(n,n)$-setoid will be called simply {\em $n$-setoid.} An $(n+1,n)$-setoid is called an {\em $n$-classoid.}
}
\end{definition}

As a justification for the term {\em classoid,} we note that there is a "replacement scheme":
if $f:A \rightarrow B$ is an extensional function from $A$, an $m$-setoid,  to  $B$ an $m$-classoid,  the image
${\rm Im}(f)$ is an $m$-setoid. This is analogous to the replacement scheme in set theory.

\medskip
\begin{example} {\em 1.  If $A \in U_n$, then $\overline{A}=(A, {\rm Id}_A(\cdot, \cdot))$ is an $n$-setoid.

2. The pair $\Omega_n = (U_n,\leftrightarrow)$, where $\leftrightarrow$ is logical equivalence, is an $n$-classoid.

3. Aczel's standard model $V$ of CZF, is built on the W-type over a universe $U_0$ with $|V| = W(U_0,T_0)$ and the equality $ =_V$ defined by bisimulation as function $|V| \rightarrow |V| \rightarrow U_0$. Thus $V=(|V|, =_V)$   forms a 0-classoid. Similarly constructing $V_k$ from a universe $U_k$, $T_k$ yields a $k$-classoid.

4. If $A$ is an $(m,n)$-setoid and $B$ is an $(m',n')$-setoid,
then the exponential $[A \rightarrow B]$ is an $(\max(m,m',n,n'),\max(m,n'))$-setoid.
In particular, if $A$ and $B$ are both $(m,n)$-setoid, then $[A \to B]$ is an
$(\max(m,n),\max(m,n))$-setoid. Thus $(m,n)$-setoids are also closed under exponents.

5. For an $n$-setoid $A$, the setoid of extensional propositional
  functions of level $n$
$$P_n(A) = [A \rightarrow \Omega_n]$$ is an $n$-classoid.

}
\end{example}

Subsetoids maybe defined following Bishop cf. \cite{palmgren2005}.

\begin{definition} {\em Let  $A$ be a setoid. A {\em subsetoid of $A$} is a setoid $\delta S$ together with an injective setoid map $\iota_S : \delta S \rightarrow A$.  An element $a : A$ is said to be a {\em member} of the subsetoid $S= (\delta S, \iota_S)$ if there is an $s:\delta S$ such that $a =_A \iota_S(s)$.  We then write 
$a \in S$ or $a \in_A S$. (Note that $s$ is unique.)  If $S$ and $T$ are subsetoids of $A$, then
we define
\begin{equation} \label{subsdef}
S \subseteq_A T \Longleftrightarrow_{\rm def} (\forall x: A)[x \in_A  S \Longrightarrow
 y \in_B T].
\end{equation}
When $A$ is clear from the context we drop this subscript.
Define also $$S \equiv_A T \Longleftrightarrow_{\rm def}  S \subseteq_A T  \land  T \subseteq_A S.$$
}
\end{definition}

\medskip 
Using the axiom of unique choice it can be seen that
\begin{equation} \label{subsdef2}
S \subseteq_A T \Longleftrightarrow (\exists f: [\delta S \rightarrow \delta T])[\iota_T \circ f = \iota_S].
\end{equation}
Here $f$ is in fact injective and unique.  
Whenever $(\delta S,\iota_S) \equiv_A (\delta T,\iota_T)$ there is a unique isomorphism $\phi : \delta S \to \delta T$ such that $\iota_T \circ \phi = \iota_S$.

\begin{definition} {\em
 Let $A$ be a fixed $(m,n)$-setoid. Define ${\rm Sub}^{m,n}_{k,\ell}(A)$ to be the type of all subsetoids $(S,\iota_S)$  such that $S$ is a $(k,\ell)$-setoid,  and the type is equipped with the equivalence relation $\equiv_A$.   Each such subsetoid is given by the data
\begin{itemize}
\item $|S| \in U_k$,  
\item $=_S : |S| \rightarrow |S| \rightarrow U_\ell$, 
\item $|\iota_S| : |S| \rightarrow |A|$
\end{itemize}
such that $(\forall x,  y : |S|)[x=_S y \leftrightarrow \iota_S(x) =_A  \iota_S(y)].$
}
\end{definition}

In this paper we will be using only the cases ${\rm Sub}^{m,m}_{m,m}(A)$ and   ${\rm Sub}^{m+1,m}_{m,m}(A)$, i.e. when $A$ is an $m$-setoid or an $m$-classoid, and we are collecting the $m$-subsetoids. However 
the levels for the general cases can be analysed as follows.

\begin{remark} {\em
The data of the subsetoid are captured by a $\Sigma$-construction in the universe of level $\max(k+1,\ell+1, m, n)$.
Two subsetoids $(S,\iota_S)$ and $(T,\iota_T)$ are equal, $(S,\iota_S) \equiv_A (T,\iota_T)$, iff  
$(\exists f: [S \rightarrow T])[\iota_T \circ f = \iota_S]$
and
$(\exists g: [T \rightarrow S])[\iota_S \circ g = \iota_T].$
Now $[S \rightarrow T]$  has type level  $\max(k,\ell)$ and $\iota_T \circ f = \iota_S$ has level $\max(k,n)$. The level of the equivalence relation is thus $\max(k, \ell, n)$.

Thus  ${\rm Sub}^{m,n}_{k,\ell}(A)$ forms a $(\max(k+1,\ell+1, m, n), \max(k, \ell, n))$-setoid.

In particular ${\rm Sub}^{m,m}_{m,m}(A)$ forms a $(m+1, m)$-setoid, i.e. an $m$-classoid. 

Further ${\rm Sub}^{m+1,m}_{m,m}(A)$ forms a $(m+1, m)$-setoid, i.e. it is also an $m$-classoid.

For $A$ an $m$-setoid or an $m$-classoid, we write ${\rm Sub}(A)$ for  ${\rm Sub}^{m,m}_{m,m}(A)$ and ${\rm Sub}^{m+1,m}_{m,m}(A)$ respectively. Thus in either case ${\rm Sub}(A)$ is an $m$-classoid.
}
\end{remark}

\begin{example} {\em  ${\rm Sub}({\mathbb V})$ is a $0$-classoid.

}
\end{example}

\subsection{Families of setoids}

\medskip
\begin{definition} {\em  Let $A$ be a setoid. A {\em proof-irrelevant setoid-family} consists of 
a family $F(a)$ of setoids indexed by $a : A$, with extensional transport 
functions $F(p) : F(a) \to F(b)$ for each proof $p : a =_A b$, that are satisfying
\begin{itemize}
\item $F(p) =_{\rm ext} F(q)$ for each pair of proofs $p, q : a =_A b$ (proof-irrelevance)
\item $F({\rm r}_a) = {\rm id}_{F(a)}$ where ${\rm r}_a : a=_A a$ is the standard proof of reflexivity. 
\item $F(p \odot q) = F(p) \circ F(q)$ if $q : a =_A b$ and $p : b =_A c$, and where 
$p \odot q : a =_A c$, using the standard proof $\odot$ of transitivity.

\end{itemize} 
}
\end{definition}

We also write $p^{-1}: b =_A a$ for $a =_A b$ using the standard proof $(-)^{-1}$ of symmetry.
Note that by functoriality and proof-irrelevance
 $$F(p) \circ F(p^{-1}) = F(p^{-1}) \circ F(p) = {\rm id}_{F(a)}$$
 so each $F(p)$ is an isomorphism.

\medskip Below we will refer to a proof-irrelevant family of setoids as just a {\em family of setoids,} when there is no chance of confusion.

\medskip
\begin{example} \label{exkappa}
{\em The operation $\kappa$ of (\ref{can-std}) extends to a family of setoids over the classoid ${\mathbb V}$, for 
$p : \alpha =_V \beta$ we let
$\kappa(p): [\kappa(\alpha) \rightarrow \kappa(\beta)]$ be given by
$$\kappa(p)(x)= \pi_1(p_1(x))$$
where 
$$p= (p_1, p_2): ((\Pi x: A)(\Sigma y:B) f(x) =_V g(y)) \times ((\Pi y: B)(\Sigma x:A) f(x) =_V g(y))$$
assuming $\alpha= {\rm sup}(A,f)$ and $\beta= {\rm sup}(B,g)$. 
}
\end{example}

\medskip
Functions into power setoids give families of setoids as can be expected:

\medskip
\begin{example} \label{fam-from-sub} {\em Let $A$ and $X$ be setoids.
Let $F: A \rightarrow {\rm Sub}(X)$ be an extensional function. Then $F(x) = (\delta(F(x)), \iota_{F(x)})$, with $\iota_{F(x)} : \delta F(x) \to X$ injective, and for $p: x =_A y$, there
is a unique isomorphism $\phi_p : \delta(F(x)) \rightarrow  \delta(F(y))$ such
that 
\begin{equation} \label{embeq}
 \iota_{F(x)}   = \iota_{F(y)} \phi_p.
\end{equation}
Thus we obtain a proof-irrelevant family $F^*$ of setoids over $A$ by letting:
$$F^*(x) := \delta(F(x)) \qquad F^*(p) := \phi_p.$$
The conditions of the transport function are easy to check.

}
\end{example}

For a setoid map $f: A \to X$ we say $f$  is  a  {\em global member of $F$} if for all $x : A$‚
$f(x) \in_X F(x)$. Thus for every $x:A$, there is a unique $s : \delta(F(x))$ such that
$$f(x) =_X \iota_{F(x)}(s).$$
Thus there is a unique function $f^* : (\Pi x : |A|)|F^*(x)|$ such that
$$f(x) =_X \iota_{F(x)}(f^*(x)).$$
If $p: x=_A y$, then $f(x) =_X f(y)$, so indeed
$$\iota_{F(x)}(f^*(x)) =_X \iota_{F(y)}(f^*(y)).$$
By (\ref{embeq}) we get
$$ \iota_{F(y)}(\phi_p(f^*(x))) =_X \iota_{F(y)}(f^*(y)).$$
Now since $\iota_{F(y)}$ is injective we have
$ \phi_p(f^*(x)) =_{F^*(y)} f^*(y)$, and thus
$$F^*(p)(f^*(x)) =_{F^*(y)} f^*(y).$$

This leads to the following definition:

\begin{definition} {\em
Let $G$ be a setoid family over $A$. A {\em global element of $G$} is a family $g(x) : G(x)$ of elements indexed by $x:A$, which is extensional in the sense that
$$(\forall p : x=_A y)[G(p)(g(x)) =_{G(y)} g(y)]$$
}
\end{definition}

Note that if $g=(|g|, {\rm ext}_g) : A \to B$ is an extensional function, and $F$ is a family on $B$, the we can form a family by composition $F \circ g$ on $A$ by defining
\begin{itemize}
\item  $(F\circ g)(x) := F(|g|(x))$ for $x : A$
\item  $(F\circ g)(p) := F({\rm ext}_g(p))$ for $p : x =_A y$ and $x, y : A$
\end{itemize}

If $f$ is a global element of $F$, then $f\circ g$ is a global element of $F\circ g$.

\begin{definition} {\em
For  $F$ a family on $A$, we can form the {\em dependent sum} $\Sigma(A,F)$ and the {\em dependent product setoid} $\Pi(A,F)$ as follows

$\Sigma(A,F) = ((\Sigma x: |A|)|F(x)|, \sim)$ where 
$$(x,y) \sim (u,v)  :=  (\exists p : x =_A u)[ F(p)(y) =_{B(u)} v ]$$

$\Pi(A,F) = (P, \sim)$ where 
\begin{equation} \label{pisetoid}
P :=  (\Sigma f : (\Pi x: |A|)|F(x)|)(\forall x, y : A)(\forall p : x=_A y)[F(p)(f(x)) =_{B(y)} f(y) ]
\end{equation}
$$ (f, e) \sim (g, e') := (\forall x: A)[f(x) =_{B(x)} g(y)].$$
}
\end{definition}

Note that $\Pi(A,F)$ consists of the global elements of $F$.

\medskip
Next we introduce an auxiliary notion. For a classoid $X$ and a family $H$ of setoids over $X$ we define a classoid of {\em parameterizations}
$${\rm Par}(X,H) =((\Sigma I : | X |)|[H(I) \rightarrow X]|, =_{{\rm Par}(X,H)}).$$
where the equivalence relation $ (I, f) =_{{\rm Par}(X,H)} (I',f')$ is defined as
$$(\exists p: I =_X I')(\forall x : H(I))f(x) =_X f'(H(p)(x))$$
This construction is used in (\ref{sigVV}) below.

\section{Basic types}  \label{ml0}

From the type universe Set in Agda we construct Aczel's type of iterative sets $V$
$$\infer[(V\; {\rm intro})]{{\rm sup}(A,f) : V}{ A : {\rm Set} & f: A \rightarrow V}$$
which corresponds to the Agda recursive data type definition
\begin{verbatim}
data V  : Set1 where
   sup : (A : Set) -> (f : A -> V) -> V
\end{verbatim}
(Note that $V$ lives in the next type universe Set1 of Agda (Set is Set0).)
This definition introduces two constants,  ${\tt V}$ a code for a type, and  ${\tt sup}$ an introduction constant.
Explained in terms of LF (Section \ref{LFAgda}) one can say that Agda automatically generates
the corresponding elimination constant and the associated the computational equality.

Define operations to extract the {\em index type} $A$ and the {\em $a$th element $f(a)$} from the set ${\rm sup}(A,f)$:
$$\# {\rm sup}(A,f) = A \qquad {\rm sup}(A,f) \blacktriangleright a = f(a).$$
The familiar set-theoretic construction  $< a, b > \; = \{\{a\},\{a, b\}\}$ of ordered pairs is used.

The natural numbers in $V$ may be constructed as the set
$${\rm natV}= {\rm sup}(N, {\rm nV})$$
where ${\rm nV}(0) = \emptyset$, and ${\rm nV}(s(m)) = \{{\rm nV}(m)\}$.
Then it can readily be shown that $\kappa({\rm natV})$ is isomorphic to the standard setoid ${\mathbb N}$ of natural numbers.

The set-theoretic version of the $\Sigma$-construction is, for $a: V$, $g : [\kappa(a) \rightarrow {\mathbb V}]$,
$${\rm sigmaV}(a,g) = {\rm sup}((\Sigma y: \#(a)) \#(g(y)), \lambda u. <a ‣ (\pi_1(u)) , (g (\pi_1(u))) ‣ (\pi_2(u)) >)$$
or in expressed in Agda code:

\begin{verbatim}
sigmaV : (a : V) -> (g : setoidmap1 (κ a) VV) -> V
sigmaV a g =  
       sup (Σ (# a) (\y -> # (g • y))) 
              (\u -> < a ‣ (pj1 u) , (g • (pj1 u)) ‣ (pj2 u) >)
\end{verbatim}
Here ${\tt pj1\; u}$ and ${\tt pj2\; u}$ denote the first and second projection of the $\Sigma$-type, respectively.
Further ${\tt VV}$ is the classoid ${\mathbb V}$, and ${\tt setoidmap1}$ is type of extensional maps. The operator $\cdot$ indicates application of such maps.

Viewing $\kappa$ as a family of setoids over ${\mathbb V}$ we define
\begin{equation} \label{sigVV}
{\rm sigma}{\mathbb V} = \lambda u.{\rm sigmaV}(\pi_1(u), \pi_2(u))  : [{\rm Par}({\mathbb V},\kappa) \rightarrow {\mathbb V}].
\end{equation}

The set-theoretic $\Pi$-construction is more involved. For $a : V$ and $g: [\kappa(a) \rightarrow {\mathbb V}]$ define 
\begin{eqnarray*}
\text{piV-iV}(a, g) &=& (\Sigma f : (\Pi x : \#(a))\#(g(a))) \\
& & \quad (\forall x, y : \#(a))(\forall p : x =_{\kappa(a)} y) (\kappa \circ g)(p)(f(x)) =_{(\kappa \circ g)(y)} f(y) \\
\text{piV-bV}(a, g) &=& \lambda h. {\rm sup}(\#(a),  (\lambda x.  < a ‣ x , g(x) ‣ (\pi_1(h)(x)) >)) \\
\text{piV}(a, g) &= & {\rm sup}(\text{piV-iV}, \text{piV-bV}(a, g)) \\
\end{eqnarray*}
The first type $\text{piV-iV}(a, g)$ singles out the extensional functions employing a $\Sigma$-type just as in (\ref{pisetoid}). The branching function 
$\text{piV-bV}$ then transforms such an extensional function to its graph in terms of set-theoretic pairs.
Similarly to the sigma-construction we define:
$${\rm pi}{\mathbb V} = \lambda u.{\rm piV}(\pi_1(u), \pi_2(u))  : [{\rm Par}({\mathbb V},\kappa) \rightarrow {\mathbb V}]$$

\begin{remark} {\em The actual  formalization uses the built-in $\Pi$-type of Agda {\tt (x:A) -> B} which may in contrast to the standard $\Pi$-type of type theory satisfy the $\eta$-rule. The Agda code is: 
\begin{verbatim}
piV-iV : (a : V) -> (g : setoidmap1 (κ a) VV) -> Set
piV-iV a g = 
      Σ ((x : # a) ->  # (g • x))
         (\f -> (x y : # a) ->
             (p : < κ a > x ~ y) ->
              < (κ° g)  § y > (ap (κ° g ± p) (f x)) ~ f y)
\end{verbatim}
}
\end{remark}

The interpretation of the extensional identity is as expected very simple: for $a:V$ and $x, y : \kappa(a)$, let
$${\rm idV}(a, x, y) = {\rm sup}((a ‣ x  =_V a ‣ y), (\lambda u.a ‣ x)).$$

\section{Universes}  \label{mlomega}

We use the type universe Set as a superuniverse \cite{palmgren1998}.  Agda's {\rm data} construct allows building universes via a so-called {\em simultaneous inductive recursive definition} \cite{dybjer2000}, such a definition has two parts, one inductive part which builds up the data part (${\tt Uo}$ below), and a second part which defines a function (${\tt To}$ below) recursion on the data part. These parts may depend mutually on each other, as in the example below, where it is crucial.

\begin{verbatim}
mutual
  data Uo (A : Set) (B : A -> Set) : Set where
     n₀ : Uo A B
     n₁ : Uo A B
     n : Uo A B
     ix : Uo A B
     lft : A -> Uo A B
     _⊕_ : Uo A B -> Uo A B -> Uo A B  
     _⊗_ : Uo A B -> Uo A B -> Uo A B  
     σ : (a : Uo A B) -> (To a -> Uo A B) -> Uo A B
     π : (a : Uo A B) -> (To a -> Uo A B) -> Uo A B
     w : (a : Uo A B) -> (To a -> Uo A B) -> Uo A B

  To : {A : Set} {B : A -> Set} ->  Uo A B -> Set
  To n₀              = N₀
  To n₁              = N₁
  To n               = N
  To {A} {B} ix      = A
  To {A} {B} (lft a)  = B a
  To (a ⊕ b)         = To a + To b
  To (a ⊗ b)         = prod (To a) (To b)
  To (σ a b)         = Σ (To a) (\x -> To (b x))
  To (π a b)         = (x  : To a) -> To (b x)  
  To (w a b)         = W (To a) (\x -> To (b x))
\end{verbatim}

To explain  the above, we note that the universe
$$\infer{{\rm Uo}(A,(x)B)}{ A \; {\rm type} & x:A \Longrightarrow B \; {\rm type}} \qquad
\infer{{\rm To}(A,(x)B,a) \; {\rm type}}{a : A}
$$
has the same closure rules as type universes à la Tarski in \cite{martinlof1984}. In addition  it has constructors
for lifting a given family $A,(x)B$ into the universe
$$\infer{{\rm ix} : {\rm Uo}(A,(x)B)}{} \qquad \infer{{\rm To}(A,(x)B,{\rm ix})=A}{}
$$
$$\infer{{\rm lft}(a) : {\rm Uo}(A,(x)B)}{ a : A} \qquad \infer{{\rm To}(A,(x)B,{\rm lft}(a))=B(a/x)}{a : A}
$$
See \cite{palmgren1998} for details.

Considering that the set universe $V$ can be obtained by applying a W-type
\begin{verbatim}
data W (A : Set) (B : A -> Set) : Set where
   sup : (a : A) -> (b : B a -> W A B) -> W A B
\end{verbatim}
to a type universe \cite{martinlof1984, aczel1982}, we get a method for constructing a hierarchy of Aczel universes. This gives us a set universe ${\rm sV}(I, F)$ for each family of types $I, F$.
\begin{verbatim}
sV : (I : Set) -> (F : I -> Set)  -> Set
sV I F = W (Uo I F) (To {I} {F})
\end{verbatim}
The elements of the small set universe ${\rm sV}(I, F)$ can be embedded into $V$
\begin{verbatim}
emb :   (I : Set) -> (F : I -> Set) -> sV I F -> V 
emb I F (sup A f) = sup (To {I} {F} A) (\x -> emb I F (f x))
\end{verbatim}
and  they form a set ${\rm uV}(I,F)$ in $V$
\begin{verbatim}
uV : (I : Set) -> (F : I -> Set) -> V
uV I F = sup (sV I F) (emb I F)
\end{verbatim}
We can think of ${\rm uV}(I,F)$ as a constructive version of an inaccessible \cite{rathjen+griffor+palmgren}.
Now, iterating the universe building operator
\begin{verbatim}
mutual

  I- : (k : N) -> Set
  I- O = I0
  I- (s k) = Uo (I- k) (F- k)

  F- : (k : N) -> I- k -> Set
  F- O = F0
  F- (s k) = To {I- k} {F- k}
\end{verbatim}
(here $I_0$, $F_0$ is an empty family) we then obtain an infinite hierarchy of inaccessibles  
$$V_k = {\rm uV}(I_k, F_k)$$ in $V$ such that $V_k \in V_{k+1}$. Each is  a transitive set so 
$V_k \subseteq V_{k+1} \subseteq V$. This will be the basis for the interpretation of the hierarchy of universes in of extensional type theory \cite{martinlof1984}.

\section{Bracket and quotient types} \label{brack}

The {\em bracket type} is a type construction  which to any type $A$ introduces a type $[A]$ whose elements are all definitionally equal (Awodey and Bauer \cite{awodey+bauer}). The idea is that $[A]$ is the proposition corresponding to $A$, and $[A]$ is inhabited if and only if $A$ is inhabited, but $[A]$ does not distinguish the proof objects. These properties are expressed by introduction and elimination rules, and some further equalities. See Section \ref{bracket} (where the notation  ${\rm Br}(A)$ is used for $[A]$).

A corresponding set-theoretic construction we use for the interpretation is the "set squasher". If $\alpha= {\rm sup}(A,f)$ is an arbitrary set, then its squashed version is
$${\rm Sq}(\alpha) =_{\rm def} {\rm sup}(A,\lambda x. \emptyset).$$
Clearly all its elements must be equal (to $\emptyset$), and also ${\rm Sq}(\alpha)$ has an element just in case $\alpha$ has an element.

Bracket types are one extreme form of quotient types. In fact CZF and hence its models admit general quotient sets, see e.g. \cite{aczel+rathjen}. Quotient rules for extensional type theory have been formulated by  Hofmann \cite[Ch.\ 5.1.5]{hofmann1997a}, Maietti \cite{maietti2005} and for HoTT in \cite{hott}.

\section{Interpretation} \label{interpretation}

Now we fix the interpretation. Define the judgements on the left to have the meaning of those on the right.
\begin{equation} \label{newinterp2}
\begin{array}{llll}
&\Gamma  \; {\rm context}  &\qquad &   \Gamma : {\mathbb V} \\
&\Gamma \Longrightarrow A \; {\rm type}  &\qquad &   A : [\kappa(\Gamma) \rightarrow {\mathbb V}]\\
&\Gamma \Longrightarrow A  == B  &&   A =_{\rm ext} B : [\kappa(\Gamma) \rightarrow  {\mathbb V}]\\
&\Gamma \Longrightarrow a  \; {\rm raw} && a : [\kappa(\Gamma) \rightarrow {\mathbb V}] \\
&\Gamma \Longrightarrow a  ::  A &&  \forall x : \kappa(\Gamma), a(x) \in_V A(x) \\
& \Gamma \Longrightarrow a == b :: A &&   \forall x : \kappa(\Gamma), a(x)=_V b(x) \in_V A(x)
\end{array}
\end{equation}
Those on the right are judgements in Agda about setoids. As usual we assume that judgements satisfy all their presuppositions.

Further we introduce judgements for substitutions between contexts, and their corresponding interpretations
\begin{equation} \label{newinterp3}
\begin{array}{llll}
&f: \Delta \longrightarrow \Gamma   &\qquad &   f : [\kappa(\Delta) \rightarrow \kappa(\Gamma)] \\
&f == g: \Delta \longrightarrow \Gamma   &\qquad &   f  =_{\rm ext} g : [\kappa(\Delta) \rightarrow \kappa(\Gamma)] \\
\end{array}
\end{equation}
The interpretation of application of substutitions to (raw) types and terms is given by composition
\begin{equation} \label{newinterp4}
\begin{array}{llll}
&a[f]  &\qquad &   a \circ f  \\
&A[f]   &\qquad &   A \circ f  \\
\end{array}
\end{equation}
Composition of substitutions is interpreted as composition of maps
\begin{equation} \label{newinterp5}
\begin{array}{llll}
&f \frown g  &\qquad &   f \circ g 
\end{array}
\end{equation}

Next the operations for context extension  $\rhd$ and the display map/left projection $\downarrow$, last variable/right projection ${\rm v}$, and extension of substitutions $ \langle \; , \; \rangle$ are defined.
\begin{equation} \label{newinterp6}
\begin{array}{llll}
&\Gamma \rhd A  &\qquad &   {\rm sigmaV}(\Gamma, A) \\
& \downarrow(A) : \Gamma \rhd A \longrightarrow \Gamma & \qquad & \pi_1 : [\kappa({\rm sigmaV}(\Gamma, A)) \rightarrow \kappa(\Gamma)] \\
& \Gamma \rhd A \Longrightarrow {\rm v}_A \; {\rm raw} & \qquad & \pi_2 : [\kappa({\rm sigmaV}(\Gamma, A)) \rightarrow {\mathbb V}] \\
& \langle f , a \rangle_{A,p} : \Delta \longrightarrow \Gamma \rhd A& \qquad & 
\lambda u. < f(u)  , \pi_1(p(u)) >\; : [\kappa(\Delta) \rightarrow \kappa({\rm sigmaV}(\Gamma, A)) ] 
 \\
\end{array}
\end{equation}
here $p : (\Delta \Longrightarrow a :: A[f])$

Finally we may introduce a judgement for equality of contexts which is interpreted as equality of sets
\begin{equation} \label{newinterp7}
\begin{array}{llll}
&\Delta == \Gamma   &\qquad &   \Delta =_V \Gamma \\
\end{array}
\end{equation}
Now by Example \ref{exkappa} each $p :  \Delta =_V \Gamma$ gives an isomorphism 
$$\phi_p =_{\rm def} \kappa(p) : [\kappa(\Delta) \rightarrow \kappa(\Gamma)]$$
 which is independent of $p$ and functorial in $p$.
Moreover it has the property that
$$\infer{\Gamma \rhd A == \Delta \rhd B}{p :  (\Delta == \Gamma) & \Gamma \Longrightarrow  A \; {\rm type}& \Delta \Longrightarrow B  \; {\rm type} &  \Gamma \Longrightarrow A == B [\phi_p]}.$$

\medskip
Some remarks about the notation to guide reading of the code.
The interpretation will mainly use 0-setoids and 0-classoids, simply called setoids and classoids. Due to some limitations of Agda notation (no subscripts) we use the following notation for $a =_A a'$ and $b=_B b'$, when $A$ and $B$ are respectively setoids and classoids
\begin{verbatim}
     < A > a ~ a'        < < B > > b ~ b'
\end{verbatim}
The underlying types $|A|$ and $|B|$ are denoted respectively
\begin{verbatim}
     || A ||           ||| B |||
\end{verbatim}
When $A, A'$ are setoids, and $B,B'$ are classoids, we use the following notations 
for $|[A,A']|$,  $|[A,B]|$ and $|[B,B']|$
\begin{verbatim}
     setoidmap A A'      setoidmap1 A B     setoidmap11 A B
\end{verbatim}

\section{Interpreted rules} \label{listrules}

The following is a list of the interpreted  rules of the formalization (Section \ref{overview}). The rule names refer to the Agda code.

We recall that the judgement forms are
$$
\begin{array}{ll}
\Gamma \; \text{context} & \Gamma \Longrightarrow A \; \text{type}\\
\Gamma == \Delta &  \Gamma \Longrightarrow A == B\\
f: \Gamma \longrightarrow \Delta & \Gamma \Longrightarrow a :: A\\
f == g : \Gamma \longrightarrow \Delta &  \Gamma \Longrightarrow a==b :: A
\end{array}
$$

The following presupposition rules are valid in the model
$$\infer{\Gamma \; \text{context}}{\Gamma== \Delta} 
\qquad \infer{\Delta \; \text{context}}{\Gamma== \Delta} 
\qquad \infer{\Gamma \; \text{context}}{f: \Gamma \longrightarrow \Delta} 
\qquad \infer{\Delta \; \text{context}}{f: \Gamma \longrightarrow \Delta} 
\qquad \infer{\Gamma \; \text{context}}{\Gamma \Longrightarrow A \; \text{type}}$$

$$\infer{\Gamma \Longrightarrow A \; \text{type}}{\Gamma \Longrightarrow A == B}
\qquad\infer{\Gamma \Longrightarrow B \; \text{type}}{\Gamma \Longrightarrow A == B}
\qquad \infer{\Gamma \Longrightarrow A \; \text{type}}{\Gamma \Longrightarrow a :: A}
$$

$$\infer{f : \Gamma \longrightarrow \Delta}{f == g : \Gamma \longrightarrow \Delta}
\qquad \infer{g : \Gamma \longrightarrow \Delta}{f == g : \Gamma \longrightarrow \Delta}
\qquad \infer{\Gamma \Longrightarrow a :: A}{\Gamma \Longrightarrow a==b :: A}
\qquad \infer{\Gamma \Longrightarrow b :: A}{\Gamma \Longrightarrow a==b :: A}
$$

\subsection{Substitutions and general equality rules}

$$\infer{\Gamma == \Gamma}{\Gamma \; {\rm context}} 
\qquad \infer{\Delta == \Gamma}{\Gamma == \Delta}
\qquad \infer{\Gamma == \Phi}{\Gamma == \Delta & \Delta == \Phi} $$

$$\infer{f== f : \Gamma \longrightarrow \Delta}{f:\Gamma \longrightarrow \Delta} 
\qquad \infer{g== f : \Gamma \longrightarrow \Delta}{f== g:\Gamma \longrightarrow \Delta}
\qquad \infer{f== h : \Gamma \longrightarrow \Delta}{f== g:\Gamma \longrightarrow \Delta & g== h:\Gamma \longrightarrow \Delta}$$

$$\infer{{\rm id}_\Gamma : \Gamma \longrightarrow \Gamma}{\Gamma \; {\rm context}} 
\qquad \infer{f \frown g: \Gamma \longrightarrow \Phi}{g : \Gamma \longrightarrow \Delta &  f: \Delta \longrightarrow  \Phi} $$

$$ \infer{g \frown {\rm id}_\Gamma == g: \Gamma \longrightarrow \Delta}{g : \Gamma \longrightarrow \Delta} 
\qquad 
\infer{{\rm id}_\Delta \frown g == g: \Gamma \longrightarrow \Delta}{g : \Gamma \longrightarrow \Delta}$$

$$\infer{(f \frown g) \frown h == f \frown (g \frown h) : \Gamma \longrightarrow \Xi }{h : \Gamma \longrightarrow \Delta & g : \Delta \longrightarrow \Phi &  f: \Phi \longrightarrow  \Xi}$$

$$\infer{f \frown g == f' \frown g': \Gamma \longrightarrow \Phi}{g == g' : \Gamma \longrightarrow \Delta &  f == f': \Delta \longrightarrow  \Phi}$$

$$\infer[(\text{subst-trp})]{\phi_p : \Gamma \longrightarrow \Delta}{p : \Gamma == \Delta}
\qquad \infer[(\text{subst-trp-irr})]{\phi_p == \phi_q : \Gamma \longrightarrow \Delta}{p : \Gamma == \Delta & q : \Gamma == \Delta}$$

$$\infer[(\text{subst-trp-id})]{\phi_p = {\rm id}_\Gamma: \Gamma \longrightarrow \Gamma}{p : \Gamma == \Gamma}
\qquad \infer[(\text{subst-trp-fun})]{\phi_q \frown \phi_p == \phi_r : \Gamma \longrightarrow \Phi}{p : \Gamma == \Delta & q : \Delta == \Phi & r : \Gamma == \Phi}$$

$$\infer[(\text{tyrefl})]{\Gamma \Longrightarrow A == A}{\Gamma \Longrightarrow A \; {\rm type}} 
\qquad \infer[(\text{tysym})]{\Gamma \Longrightarrow B == A}{\Gamma \Longrightarrow A == B}
\qquad  \infer[(\text{tytra})]{\Gamma \Longrightarrow A == C}{\Gamma \Longrightarrow A == B & \Gamma \Longrightarrow B == C} $$

$$\infer[(\text{tmrefl})]{\Gamma \Longrightarrow a == a :: A}{\Gamma \Longrightarrow a :: A } 
\qquad \infer[(\text{tmsym})]{\Gamma \Longrightarrow b == a :: A}{\Gamma \Longrightarrow a == b :: A}$$ 

$$\infer[(\text{tmtra})]{\Gamma \Longrightarrow a == c :: A}{\Gamma \Longrightarrow a == b :: A & \Gamma \Longrightarrow b == c :: A} $$

$$\infer[(\text{elttyeq})]{\Gamma \Longrightarrow a :: B }{\Gamma \Longrightarrow a :: A & \Gamma \Longrightarrow A == B } 
\qquad \infer[(\text{elteqtyeq})]{\Gamma \Longrightarrow a == b :: B }{\Gamma \Longrightarrow a == b:: A & \Gamma \Longrightarrow A == B } $$

$$\infer{\Delta \Longrightarrow A[f] \; {\rm type}}{\Gamma \Longrightarrow A \; {\rm type} & f : \Delta \longrightarrow \Gamma}$$

$$\infer[(\text{tyeq-subst})]{\Delta \Longrightarrow A[f] == B[f]}{\Gamma \Longrightarrow A == B & f : \Delta \longrightarrow \Gamma}
\qquad \infer[(\text{tyeq-subst2})]{\Delta \Longrightarrow A[f] == A[g]}{\Gamma \Longrightarrow A \; {\rm type} & f == g : \Delta \longrightarrow \Gamma}$$

$$\infer[(\text{tysubst-id})]{\Gamma \Longrightarrow A[{\rm id}_\Gamma] == A}{\Gamma \Longrightarrow A\; {\rm type}}
\qquad \infer[(\text{tysubst-com})]{\Phi \Longrightarrow A[f \frown g] == A[f][g]}{\Gamma \Longrightarrow A \; {\rm type}& g : \Phi \longrightarrow \Delta & f : \Delta \longrightarrow \Gamma}$$

$$\infer[(\text{elt-subst})]{\Delta \Longrightarrow a[f] :: A[f]}{\Gamma \Longrightarrow a :: A & f : \Delta \longrightarrow \Gamma}$$

$$\infer[(\text{elteq-subst})]{\Delta \Longrightarrow a[f] == b[f] :: A[f]}{\Gamma \Longrightarrow a == b :: A & f : \Delta \longrightarrow \Gamma}
\qquad \infer[(\text{elteq-subst2})]{\Delta \Longrightarrow a[f] == a[g] :: A[f]}{\Gamma \Longrightarrow a :: A & f == g : \Delta \longrightarrow \Gamma}$$

$$\infer[(\text{eltsubst-id})]{\Gamma \Longrightarrow a[{\rm id}_\Gamma] == a :: A}{\Gamma \Longrightarrow a :: A}
\qquad \infer[(\text{eltsubst-com})]{\Phi \Longrightarrow a[f \frown g] == a[f][g] ::  A[f \frown g]}{\Gamma \Longrightarrow a :: A & g : \Phi \longrightarrow \Delta & f : \Delta \longrightarrow \Gamma}$$

\subsection{Context extension and associated rules}

$$\infer{\emptyctx \; {\rm context}}{ }\qquad \infer{\Gamma \rhd A \; {\rm context}}{\Gamma \Longrightarrow A \; {\rm type}}$$

$$\infer[(\text{ext-eq'})]{\Gamma \rhd A  ==  \Gamma \rhd B}{\Gamma \Longrightarrow A \; {\rm type} & \Gamma \Longrightarrow B \; {\rm type}  & \Gamma \Longrightarrow A == B}$$

$$\infer[(\text{ext-eq''})]{\Gamma \rhd A  ==  \Delta \rhd B}{\Gamma \Longrightarrow A \; {\rm type} & \Delta \Longrightarrow B \; {\rm type} & p : (\Gamma == \Delta) & \Gamma \Longrightarrow A == B [\phi_p]}$$

$$\infer[(\downarrow)]{\downarrow A : \Gamma \rhd A \longrightarrow \Gamma}{\Gamma \Longrightarrow A \; {\rm type}}$$

$$\infer[(\downarrow\text{cong})]{\varphi_p \frown (\downarrow A) == 
(\downarrow B) \frown \varphi(\text{ext-eq''}(A,B,p,q)) : \Gamma \rhd A \longrightarrow \Delta }{\Gamma \Longrightarrow A \; {\rm type} & \Delta \Longrightarrow B \; {\rm type} & p : (\Gamma == \Delta) &  q : (\Gamma \Longrightarrow A == B [\phi_p])}$$

$$\infer[(\text{asm})]{\Gamma \rhd A \Longrightarrow  {\rm v}_A :: A[\downarrow A]}{\Gamma \Longrightarrow A \; {\rm type}}$$

$$\infer[(\text{asm-cong})]{\Gamma \rhd A \Longrightarrow  {\rm v}_A  == {\rm v}_B [\varphi(\text{ext-eq''}(A,B,p,q))] :: A[\downarrow A]}{\Gamma \Longrightarrow A \; {\rm type} & \Delta \Longrightarrow B \; {\rm type} & p : (\Gamma == \Delta) &  q : (\Gamma \Longrightarrow A == B [\phi_p])}$$

$$\infer[(\text{ext})]{\langle f, a\rangle_p : \Delta \longrightarrow \Gamma \rhd A }{f: \Delta \longrightarrow \Gamma & \Gamma \Longrightarrow A \; {\rm type} &  p : (\Delta \Longrightarrow a :: A[f]) }$$

$$\infer[(\text{ext-irr})]{\langle f, a\rangle_p == \langle f, a\rangle_q : \Delta \longrightarrow \Gamma \rhd A }{f: \Delta \longrightarrow \Gamma & \Gamma \Longrightarrow A \; {\rm type} &  p : (\Delta \Longrightarrow a :: A[f]) &  q : (\Delta \Longrightarrow a :: A[f]) }$$

$$\infer[(\text{ext-cong})]{\langle f, a\rangle_p == \langle g, b\rangle_q : \Delta \longrightarrow \Gamma \rhd A }{\begin{array}{l} f == g: \Delta \longrightarrow \Gamma \\
 \Gamma \Longrightarrow A \; {\rm type} \\
  p : (\Delta \Longrightarrow a :: A[f]) \\  q : (\Delta \Longrightarrow b :: A[g]) \\
  r : (\Delta \Longrightarrow a == b:: A[f]) 
  \end{array} }$$

$$\infer[(\text{ext-prop1})]{(\downarrow A) \frown \langle f, a\rangle_p == f : \Delta \longrightarrow \Gamma }{f: \Delta \longrightarrow \Gamma & \Gamma \Longrightarrow A \; {\rm type} &  p : (\Delta \Longrightarrow a :: A[f]) }$$

$$\infer[(\text{ext-prop2})]{ \Delta  \Longrightarrow {\rm v}_A[\langle f, a\rangle_p] == a  :: A[f] }{f: \Delta \longrightarrow \Gamma & \Gamma \Longrightarrow A \; {\rm type} &  p : (\Delta \Longrightarrow a :: A[f]) }$$

$$\infer[(\text{ext-prop3})]{\langle \downarrow A, {\rm v}_A\rangle_p == {\rm id}_{\Gamma \rhd A}: \Gamma \rhd A \longrightarrow \Gamma \rhd A }{\Gamma \Longrightarrow A \; {\rm type} & p : (\Delta \Longrightarrow {\rm v}_A :: A[\downarrow A])}$$

$$\infer{\langle f, a\rangle_p \frown h = \langle f \frown h, a[h]\rangle_q : \Delta \longrightarrow \Gamma \rhd A }{ \begin{array}{l} h:  \Theta \longrightarrow \Delta \\ f: \Delta \longrightarrow \Gamma \\ \Gamma \Longrightarrow A \; {\rm type} \\  p : (\Delta \Longrightarrow a :: A[f]) \\  q : (\Delta \Longrightarrow a[h] :: A[f \frown h]) 
\end{array}}$$

\medskip
Two derived rules:


$$\infer{{\rm els}(p) : \Gamma \longrightarrow \Gamma \rhd A}{p: (\Gamma \Longrightarrow a :: A)}
\qquad 
\infer[(\text{els-exp})]{{\rm els}(p) ==  \langle {\rm id}_\Gamma, a \rangle_p: \Gamma \longrightarrow \Gamma \rhd A}{p: (\Gamma \Longrightarrow a :: A)}
$$


$$\infer{\uparrow(A,h) : \Delta \rhd A[h] \longrightarrow \Gamma \rhd A}{\Gamma \Longrightarrow A \; {\rm type} 
      & h :  \Delta \longrightarrow \Gamma}$$
$$ \infer[(\text{qq-exp})]{\uparrow(A,h) == \langle  h \frown (\downarrow A [h]) , {\rm v}_{A [h]} \rangle_p: \Delta \rhd A[h] \longrightarrow \Gamma \rhd A}{\Gamma \Longrightarrow A \; {\rm type} 
      & h :  \Delta \longrightarrow \Gamma & p : (\Delta \Longrightarrow {\rm v}_{A [h]} :: A[ h \frown (\downarrow A [h]) ])}
$$

\subsection{Rules for particular type constructions}

The general principle of Martin-L\"of type theory is that each type construction comes with a formation rule, a finite number of introduction rules, one elimination rule, and computation rules. There maybe additional equality rules in extended theories. Moreover each constant has a congruence rule. If the theory is based on explicit substitution (as  is the case here) there also equality rules that state that substitutions commute with constants and abstractions.

\subsubsection{$\Pi$-rules}

$$\infer[(\Pi\text{-f})]{\Gamma \Longrightarrow \Pi_{\rm f}(A,B) \; {\rm type}}{\Gamma \Longrightarrow A \; {\rm type} & \Gamma \rhd A \Longrightarrow B \; {\rm type}}$$

$$\infer[(\Pi\text{-i})]{\Gamma \Longrightarrow \lambda(A,B,b) :: \Pi_{\rm f}(A,B)}{\Gamma \Longrightarrow A \; {\rm type} & \Gamma \rhd A  \Longrightarrow B \; {\rm type} & \Gamma \rhd A \Longrightarrow b :: B}$$

$$\infer[(\Pi\text{-e})]{\Gamma \Longrightarrow \app(A,B,c,p,a,q) :: B[{\rm els}(q)]}{\Gamma \Longrightarrow A \; {\rm type} & \Gamma \rhd A  \Longrightarrow B \; {\rm type} & 
p : (\Gamma \Longrightarrow c :: \Pi_{\rm f}(A,B)) & q : (\Gamma \Longrightarrow a :: A)
}$$

{\small 
$$\infer[(\Pi\text{-beta-gen})]{\Gamma \Longrightarrow \app(A,B,  \lambda(A,B,b), r,  a, q) == b[{\rm els}(q)] :: B[{\rm els}(q)] }{\Gamma \Longrightarrow A \; {\rm type} & \Gamma \rhd A  \Longrightarrow B \; {\rm type} &  r : (\Gamma \Longrightarrow \lambda(A,B,b) :: \Pi_{\rm f}(A,B)) & q : ( \Gamma \Longrightarrow a :: A) }$$
}

{\small
$$\infer[(\Pi\text{-eta-eq-gen})]{\lambda(A,B,\app(A[\downarrow(A)], B[\uparrow(A, \downarrow(A))], c [\downarrow(A)],q_2,{\rm v}_A,q_1))== c :: \Pi_{\rm f}(A,B)}{  
\begin{array}{l}
p : (\Gamma \Longrightarrow c :: \Pi_{\rm f}(A,B)) \\
q_1 : ( \Gamma \rhd A \Longrightarrow {\rm v}_A :: A[\downarrow(A)]) \\
q_2 : ( \Gamma \rhd A \Longrightarrow  c [\downarrow(A)] :: \Pi_{\rm f}(A[\downarrow(A)],B[\uparrow(A, \downarrow(A))])) 
\end{array}
}$$
}

$$\infer[(\Pi\text{-f-sub})]{\Delta \Longrightarrow \Pi_{\rm f}(A,B) [h] == \Pi_{\rm f}(A[h],B[\uparrow(A,h)])}{\Gamma \Longrightarrow A \; {\rm type} & \Gamma \rhd A \Longrightarrow B \; {\rm type} & h: \Delta \longrightarrow \Gamma}$$

$$\infer[(\text{lambda-sub})]{\Delta \Longrightarrow \lambda(A,B,b) [h] ==  \lambda(A[h],B[\uparrow(A,h)], b[\uparrow(A,h)]):: \Pi_{\rm f}(A,B)[h]}{\Gamma \Longrightarrow A \; {\rm type} & \Gamma \rhd A \Longrightarrow B \; {\rm type}  &  \Gamma \rhd A \Longrightarrow b :: B & h: \Delta \longrightarrow \Gamma}$$

{\small
$$\infer[(\Pi\text{-e-sub-gen})]{\Gamma \Longrightarrow \app(A,B,c,p,a,q)[h] == \app(A[h],B[\uparrow(A,h)],c[h],r_1,a[h],r_2) :: B[{\rm els}(q)][h]}{ 
\begin{array}{l}
p : (\Gamma \Longrightarrow c :: \Pi_{\rm f}(A,B)) \\
 q : (\Gamma \Longrightarrow a :: A) \\
  h: \Delta \longrightarrow \Gamma \\ 
  r_1 :(\Delta \Longrightarrow c[h] :: \Pi_{\rm f}(A[h],B[\uparrow(A,h)])) \\
   r_2 :(\Delta \Longrightarrow a[h] :: A[h])
\end{array}}$$
}

$$\infer[(\Pi_{\rm f}\text{-cong})]{\Gamma \Longrightarrow \Pi_{\rm f}(A,B)  == \Pi_{\rm f}(A',B')}{
\begin{array}{l}
p: (\Gamma \Longrightarrow A == A') \\
 \Gamma \rhd A \Longrightarrow B \; {\rm type}  \\
  \Gamma \rhd A' \Longrightarrow B' \; {\rm type}  \\
   \Gamma \rhd A \Longrightarrow B == B'[\phi(\text{ext-eq'}(A , A', p))]
\end{array}}$$

$$\infer[(\Pi\text{-xi})]{\Gamma \Longrightarrow \lambda(A,B,b) == \lambda(A,B,b') :: \Pi_{\rm f}(A,B)}{\Gamma \Longrightarrow A \; {\rm type}& \Gamma \rhd A  \Longrightarrow B \; {\rm type}& \Gamma \rhd A \Longrightarrow b == b' :: B}$$

$$\infer[(\Pi\text{-e-cong})]{\Gamma \Longrightarrow \app(A,B,c,p,a,q)  == \app(A,B,c',p',a',q') :: B [ {\rm els}(q) ]}
    {\begin{array}{l}
    p : (\Gamma \Longrightarrow c ::  \Pi_{\rm f}(A,B))\\
    p' : (\Gamma \Longrightarrow c' ::  \Pi_{\rm f}(A,B))\\
     \Gamma \Longrightarrow c == c' ::  \Pi_{\rm f}(A,B)\\
    q : (\Gamma \Longrightarrow a ::  A)\\
    q' : (\Gamma \Longrightarrow a' ::  A)\\
     \Gamma \Longrightarrow a == a' ::  A\\
     \end{array}}
$$
      
\subsubsection{Id-rules}

$$\infer[(\text{ID})]{\Gamma \Longrightarrow {\rm ID}(A,a,p,b,q) \; {\rm type}}{\Gamma \Longrightarrow A \; {\rm type} & p: (\Gamma \Longrightarrow a :: A)& q: (\Gamma \Longrightarrow b :: A)}$$

$$\infer[(\text{ID-i})]{\Gamma \Longrightarrow {\rm rr}(a) :: {\rm ID}(A,a,p,a,p)}{p: (\Gamma \Longrightarrow a :: A)}$$

$$\infer[(\text{ID-e})]{\Gamma \Longrightarrow a == b :: A}{p: (\Gamma \Longrightarrow a :: A)& q: (\Gamma \Longrightarrow a :: A) & \Gamma \Longrightarrow t :: {\rm ID}(A,a,p,b,q)}$$

$$\infer[(\text{ID-uip})]{\Gamma \Longrightarrow t == {\rm rr}(a) :: {\rm ID}(A,a,p,a,q)}{p: (\Gamma \Longrightarrow a :: A)& q: (\Gamma \Longrightarrow a :: A) & \Gamma \Longrightarrow t :: {\rm ID}(A,a,p,a,q)}$$

$$\infer[(\text{ID-sub-gen})]{\Delta \Longrightarrow {\rm ID}(A,a,p_a,b,p_b)[ h ] == {\rm ID}(A[h], a[h],p,b[h],q)}{\begin{array}{l} 
h : \Delta \longrightarrow \Gamma \\
p_a: (\Gamma \Longrightarrow a :: A) \\ p_b: (\Gamma \Longrightarrow b :: A) \\
p: (\Delta \Longrightarrow a [h]  :: A[h] )\\
q: (\Delta \Longrightarrow b [h]  :: A[h]) 
\end{array}}$$

$$\infer[(\text{rr-sub})]{\Delta \Longrightarrow {\rm rr}(a) [h] == {\rm rr}(a [h]) :: {\rm ID}(A,a,p,a,p) [h]}
{ h :  \Delta \longrightarrow \Gamma &
\Gamma \Longrightarrow A \; {\rm type} &
p: (\Gamma \Longrightarrow a :: A) &
}$$

$$\infer[(\text{ID-cong})]{\Gamma \Longrightarrow {\rm ID}(A,a,p_a,b,p_b) ==   {\rm ID}(A',a',p_{a'},b',p_{b'})}{\begin{array}{l} 
p_a: (\Gamma \Longrightarrow a :: A) \\ 
p_{a'}: (\Gamma \Longrightarrow a' :: A') \\ 
p_b: (\Gamma \Longrightarrow b :: A) \\
p_{b'}: (\Gamma \Longrightarrow b' :: A') \\
\Gamma \Longrightarrow A == A' \\
\Gamma \Longrightarrow a == a' :: A \\
\Gamma \Longrightarrow b == b' :: A \\
\end{array}}$$

$$\infer[(\text{rr-cong})]{\Gamma \Longrightarrow {\rm rr}(a) == {\rm rr}(b) :: {\rm ID}(A,a,p,a,p)}
{ 
p: (\Gamma \Longrightarrow a :: A) &
\Gamma \Longrightarrow a == b :: A &
}$$

\subsubsection{$\Sigma$-rules}

$$\infer[(\Sigma\text{-f})]{\Gamma \Longrightarrow \Sigma_{\rm f}(A,B) \; {\rm type}}{\Gamma \Longrightarrow A \; {\rm type} & \Gamma \rhd A \Longrightarrow B \; {\rm type}}$$

$$\infer[(\Sigma\text{-i})]{\Gamma \Longrightarrow {\rm pr}(a,b) :: \Sigma_{\rm f}(A,B)}{p : (\Gamma \Longrightarrow  a :: A) \ & \Gamma \Longrightarrow b :: B[{\rm els}(p)] }$$

$$\infer[(\Sigma\text{-e-1})]{\Gamma \Longrightarrow    {\rm pr}_1(c,p) :: A }{ p : (\Gamma \Longrightarrow c :: \Sigma_{\rm f}(A,B))} \qquad 
\infer[(\Sigma\text{-e-2})]{\Gamma \Longrightarrow    {\rm pr}_2(c,p) :: B[{\rm els}(q)] }{p : (\Gamma \Longrightarrow c :: \Sigma_{\rm f}(A,B)) & q : (\Gamma \Longrightarrow    {\rm pr}_1(c,p) :: A)}
$$

$$\infer[(\Sigma\text{-c-1})]{\Gamma \Longrightarrow    {\rm pr}_1({\rm pr}(a,b),q) == a :: A }{p : (\Gamma \Longrightarrow  a :: A) \ & \Gamma \Longrightarrow b :: B[{\rm els}(p)] & q: (\Gamma \Longrightarrow {\rm pr}(a,b) :: \Sigma_{\rm f}(A,B)) }$$
 
$$\infer[(\Sigma\text{-c-2})]{\Gamma \Longrightarrow    {\rm pr}_2({\rm pr}(a,b),q) == b  :: B[{\rm els}(p)] }{p : (\Gamma \Longrightarrow  a :: A) \ & \Gamma \Longrightarrow b :: B[{\rm els}(p)]  & q : (\Gamma \Longrightarrow {\rm pr}(a,b) :: \Sigma_{\rm f}(A,B))}
$$

$$\infer[(\Sigma\text{-c-eta})]{\Gamma \Longrightarrow    c == {\rm pr}({\rm pr}_1(c,p),{\rm pr}_2(c,p))  ::  \Sigma_{\rm f}(A,B) }{p : (\Gamma \Longrightarrow c :: \Sigma_{\rm f}(A,B))}
$$

$$\infer[(\Sigma_{\rm f}\text{-cong})]{\Gamma \Longrightarrow \Sigma_{\rm f}(A,B)  == \Sigma_{\rm f}(A',B')}{
\begin{array}{l}
p: (\Gamma \Longrightarrow A == A') \\
 \Gamma \rhd A \Longrightarrow B \; {\rm type}  \\
  \Gamma \rhd A' \Longrightarrow B' \; {\rm type}  \\
   \Gamma \rhd A \Longrightarrow B == B'[\phi(\text{ext-eq'}(A , A', p))]
\end{array}}$$

$$\infer[(\text{pr-cong})]{\Gamma \Longrightarrow {\rm pr}(a,b) = {\rm pr}(a',b') :: \Sigma_{\rm f}(A,B) }{
p : (\Gamma \Longrightarrow  a :: A) & \Gamma \Longrightarrow  a == a' :: A & \Gamma \Longrightarrow b == b'  :: B[{\rm els}(p)] }$$

$$\infer[(\text{pr1-cong})]{\Gamma \Longrightarrow {\rm pr}_1(c,p) == {\rm pr}_1(c',p') :: A }{
p : (\Gamma \Longrightarrow c :: \Sigma_{\rm f}(A,B)) & p' : (\Gamma \Longrightarrow c :: \Sigma_{\rm f}(A,B))  &  (\Gamma \Longrightarrow c == c' :: \Sigma_{\rm f}(A,B))}$$

$$\infer[(\text{pr2-cong})]{\Gamma \Longrightarrow {\rm pr}_2(c,p) == {\rm pr}_2(c',p') :: B [{\rm els}(r)]) }{\begin{array}{l}
p : (\Gamma \Longrightarrow c :: \Sigma_{\rm f}(A,B)) \\ 
p' : (\Gamma \Longrightarrow c' :: \Sigma_{\rm f}(A,B)) \\
  (\Gamma \Longrightarrow c == c' :: \Sigma_{\rm f}(A,B)) \\
   r : (\Gamma \Longrightarrow {\rm pr}_1(c,p) :: A )
   \end{array} }$$

$$\infer[(\Sigma\text{-f-sub})]{\Delta \Longrightarrow \Sigma_{\rm f}(A,B) [h] == \Sigma_{\rm f}(A[h],B[\uparrow(A,h)])}{\Gamma \Longrightarrow A \; {\rm type} & \Gamma \rhd A \Longrightarrow B \; {\rm type} & h: \Delta \longrightarrow \Gamma}$$

$$\infer[(\text{pr-sub})]{\Delta \Longrightarrow {\rm pr}(a,b) [h] ==  {\rm pr}(a[h], b[h]) :: \Sigma_{\rm f}(A,B)[h]}{ h: \Delta \longrightarrow \Gamma & p : (\Gamma \Longrightarrow  a :: A)  & \Gamma \Longrightarrow b :: B[{\rm els}(p)] }$$

$$\infer[(\text{pr1-sub})]{\Delta \Longrightarrow {\rm pr}_1(c,p)[h] == {\rm pr}_1(c[h],q) :: A[h] }{
h: \Delta \longrightarrow \Gamma & p : (\Gamma \Longrightarrow c :: \Sigma_{\rm f}(A,B)) & q : (\Delta \Longrightarrow c[h] :: \Sigma_{\rm f}(A[h],B[\uparrow(A,h)]))}$$

$$\infer[(\text{pr2-sub})]{\Gamma \Longrightarrow {\rm pr}_2(c,p)[h] == {\rm pr}_2(c[h],q) :: B[\uparrow(A,h)][{\rm els}(r)] }{
\begin{array}{l}
h: \Delta \longrightarrow \Gamma \\
p : (\Gamma \Longrightarrow c :: \Sigma_{\rm f}(A,B)) \\
q : (\Gamma \Longrightarrow c[h] :: \Sigma_{\rm f}(A[h],B[\uparrow(A,h)])) \\
r : (\Delta \Longrightarrow {\rm pr}_1(c,p)[h] :: A [h])
\end{array}}$$

\subsubsection{$N$-rules}

$$\infer{\Gamma \Longrightarrow {\rm Nat} \; {\rm type}}{\Gamma \; {\rm context}}$$

$$\infer[(\text{Nat-i-0})]{\Gamma \Longrightarrow 0 :: {\rm Nat}}{\Gamma \; {\rm context}}
\qquad \infer[(\text{Nat-i-s})]{\Gamma \Longrightarrow s(a) :: {\rm Nat}}{\Gamma \Longrightarrow a :: {\rm Nat}} $$

$$\infer[(\text{Nat-e})]{\Gamma \Longrightarrow {\rm Rec}(C, d, p, e, q, c, r) :: C[{\rm els}(r)] }{
\begin{array}{l}
\Gamma \rhd {\rm Nat} \Longrightarrow C \; {\rm type} \\
p : (\Gamma \Longrightarrow d :: C[{\rm els}(\text{Nat-i-0})]) \\
q : (\Gamma \rhd {\rm Nat} \rhd C  \Longrightarrow e :: C[\text{step-sub}(\Gamma)] [\downarrow(C) ]) \\
r : (\Gamma \Longrightarrow c :: {\rm Nat})
\end{array}}
$$
Here $\text{step-sub}(\Gamma): \Gamma \rhd {\rm Nat}\longrightarrow \Gamma \rhd {\rm Nat}$ is the straightforward substitution that applies the successor to the second argument.

$$\infer[(\text{Nat-c-0})]{\Gamma \Longrightarrow {\rm Rec}(C, d, p, e, q, 0, \text{Nat-i-0}) == d :: C[{\rm els}(\text{Nat-i-0})] }{
\begin{array}{l}
\Gamma \rhd {\rm Nat} \Longrightarrow C \; {\rm type} \\
p : (\Gamma \Longrightarrow d :: C[{\rm els}(\text{Nat-i-0})]) \\
q : (\Gamma \rhd {\rm Nat} \rhd C  \Longrightarrow e :: C[\text{step-sub}(\Gamma)] [\downarrow(C) ])
\end{array}}
$$

{\small
$$\infer[(\text{Nat-c-s})]{
\begin{array}{l}
\Gamma \Longrightarrow {\rm Rec}(C, d, p, e, q, s(a), \text{Nat-i-s}(a, r) ) \\
\qquad \qquad == e [\langle {\rm els}(r), {\rm Rec}(C, d, p, e, q, a, r) \rangle_{\text{Nat-e}(C, d, p, e, q, a, r)}]:: C[{\rm els}(\text{Nat-i-s}(a, r))] 
\end{array}}{
\begin{array}{l}
\Gamma \rhd {\rm Nat} \Longrightarrow C \; {\rm type} \\
p : (\Gamma \Longrightarrow d :: C[{\rm els}(\text{Nat-i-0})]) \\
q : (\Gamma \rhd {\rm Nat} \rhd C  \Longrightarrow e :: C[\text{step-sub}(\Gamma)] [\downarrow(C) ]) \\
r : (\Gamma \Longrightarrow a :: {\rm Nat})
\end{array}}
$$
}

$$\infer[(\text{Nat-i-s-cong})]{\Gamma \Longrightarrow s(a) == s(b) :: {\rm Nat}}{\Gamma \Longrightarrow a  == b:: {\rm Nat}} $$

$$\infer[(\text{Rec-cong})]{\Gamma \Longrightarrow {\rm Rec}(C, d, p, e, q, c, r) ==  {\rm Rec}(C', d', p', e', q', c', r') :: C[{\rm els}(r)]  }{
\begin{array}{l}
p : (\Gamma \Longrightarrow d :: C[{\rm els}(\text{Nat-i-0})]) \\
p' : (\Gamma \Longrightarrow d' :: C'[{\rm els}(\text{Nat-i-0})]) \\
q : (\Gamma \rhd {\rm Nat} \rhd C  \Longrightarrow e :: C[\text{step-sub}(\Gamma)] [\downarrow(C) ]) \\
q' : (\Gamma \rhd {\rm Nat} \rhd C'  \Longrightarrow e' :: C'[\text{step-sub}(\Gamma)] [\downarrow(C') ]) \\
r : (\Gamma \Longrightarrow c :: {\rm Nat}) \\
r' : (\Gamma \Longrightarrow c' :: {\rm Nat}) \\
t : (\Gamma \rhd {\rm Nat} \Longrightarrow C == C') \\
\Gamma \Longrightarrow d == d' ::C[{\rm els}(\text{Nat-i-0})] \\
\Gamma \rhd {\rm Nat} \rhd C  \Longrightarrow e == e'[ \phi(\text{ext-eq'}(C, C', t)) ]:: C[\text{step-sub}(\Gamma)] [\downarrow(C) ]  \\
\Gamma \Longrightarrow c == c' :: {\rm Nat}
\end{array}} $$

$$\infer[(\text{Nat-sub})]{\Delta \Longrightarrow {\rm Nat} [ h  ] ==  {\rm Nat} }{h: \Delta \longrightarrow \Gamma}$$

$$\infer[(\text{Nat-i-0-sub})]{\Delta \Longrightarrow 0[ h ] = 0 :: {\rm Nat}}{h: \Delta \longrightarrow \Gamma}
\qquad \infer[(\text{Nat-i-s-sub})]{\Delta \Longrightarrow s(a) [ h ] = s ([ a]) :: {\rm Nat}}{h: \Delta \longrightarrow \Gamma & \Gamma \Longrightarrow a :: {\rm Nat}} $$

{\small
$$\infer[(\text{Rec-sub})]{\Delta \Longrightarrow {\rm Rec}(C, d, p, e, q, c, r) [h] ==
 {\rm Rec}(C[\text{N-sub}(h)], d[h], p', e[\text{C-sub}(h, C)],  q', c[h], r') ::  C[{\rm els}(r)] [h] }{
\begin{array}{l}
 h : \Delta \longrightarrow \Gamma \\
 \Gamma \rhd {\rm Nat} \Longrightarrow C\; {\rm type} \\
p : (\Gamma \Longrightarrow d :: C[{\rm els}(\text{Nat-i-0})]) \\
p' : (\Gamma \Longrightarrow d [ h ] :: C[\text{N-sub}(h)][{\rm els}(\text{Nat-i-0})]) \\
q : (\Gamma \rhd {\rm Nat} \rhd C  \Longrightarrow e :: C[\text{step-sub}(\Gamma)] [\downarrow(C) ]) \\
q' : (\Delta \rhd {\rm Nat} \rhd C[\text{N-sub}(h)]  \Longrightarrow e  [\text{C-sub}(h, C) ] :: 
C[\text{N-sub}(h)][\text{step-sub}(\Delta)] [\downarrow(C[\text{N-sub}(h)]) ]) 
\end{array}
} $$
}

\medskip
\subsubsection{$N_0$-rules}

$$\infer[(\text{N0})]{\Gamma \Longrightarrow {\rm N}_0 \; {\rm type}}{\Gamma \; {\rm context}}$$

$$\infer[(\text{N0-e})]{\Gamma \Longrightarrow {\rm R}_0(C,c,r) :: C [{\rm els}(r)]}{\Gamma \rhd {\rm N}_0 \Longrightarrow C \; {\rm type} & r : (\Gamma \Longrightarrow c :: {\rm N}_0) }$$

$$\infer[(\text{R0-cong'})]{\Gamma \Longrightarrow {\rm R}_0(C,c,r) == {\rm R}_0(C',c',r'):: C [{\rm els}(r)]}{\begin{array}{l}
\Gamma \rhd {\rm N}_0 \Longrightarrow C == C'  \\ 
r : (\Gamma \Longrightarrow c :: {\rm N}_0) \\
r' : (\Gamma \Longrightarrow c' :: {\rm N}_0) \\
\Gamma \Longrightarrow c == c':: {\rm N}_0
\end{array} }$$

$$\infer[(\text{N0-sub})]{\Delta \Longrightarrow {\rm N}_0 [ h  ] ==  {\rm N}_0 }{h: \Delta \longrightarrow \Gamma}$$

$$\infer[(\text{R0-sub})]{\Delta \Longrightarrow {\rm R}_0(C,c,r)  [ h ] == {\rm R}_0(C [ \uparrow({\rm N}_0, h) ],c [ h ] ,r') :: C [{\rm els}(r)] [ h  ]
 }{\begin{array}{l}
h: \Delta \longrightarrow \Gamma \\ 
\Gamma \rhd {\rm N}_0 \Longrightarrow C \; {\rm type} \\ 
r : (\Gamma \Longrightarrow c :: {\rm N}_0)\\  
r' : (\Gamma \Longrightarrow c  [ h ]:: {\rm N}_0 [ h  ]) 
\end{array}}$$

\subsubsection{$+$-rules}

$$\infer[(\text{Sum})]{\Gamma \Longrightarrow {\rm Sum}(A,B) \; {\rm type}}{\Gamma \Longrightarrow A \; {\rm type} & \Gamma \Longrightarrow B \; {\rm type}}$$

$$\infer[(\text{lf-pf})]{\Gamma \Longrightarrow {\rm lf}(A,B,a,p) :: {\rm Sum}(A,B) }{\Gamma \Longrightarrow B \; {\rm type}  & p : (\Gamma \Longrightarrow a :: A)}$$

$$\infer[(\text{rg-pf})]{\Gamma \Longrightarrow {\rm rg}(A,B,b,q) :: {\rm Sum}(A,B) }{\Gamma \Longrightarrow A \; {\rm type} & q : (\Gamma \Longrightarrow b :: B)}$$

$$\infer[(\text{Sum-e})]{\Gamma \Longrightarrow \text{Sum-rec}(A,B,C,d,p,e,q,c,r) :: C[{\rm els}(r)]}{
\begin{array}{l}
  \Gamma \rhd {\rm Sum}(A,B) \Longrightarrow C \; {\rm type} \\
  p: (\Gamma \rhd {\rm Sum}(A,B) \Longrightarrow d :: C [\text{Sum-sub-lf}(A,B)] \\
  q: (\Gamma \rhd {\rm Sum}(A,B) \Longrightarrow e :: C [\text{Sum-sub-rg}(A,B)] \\
  r:  (\Gamma \Longrightarrow c:: {\rm Sum}(A,B))
\end{array}}$$

Here $\text{Sum-sub-lf}(A,B): \Gamma \rhd A \longrightarrow \Gamma \rhd {\rm Sum}(A,B)$ and
$\text{Sum-sub-rg}(A,B): \Gamma \rhd B \longrightarrow \Gamma \rhd {\rm Sum}(A,B)$ are defined from ${\rm lf}$ and ${\rm rg}$ using $\langle,\rangle$ in the straightforward way. 

$$\infer[(\text{Sum-c1})]{\Gamma \Longrightarrow \text{Sum-rec}(A,B,C,d,p,e,q,{\rm lf}(A,B,a,r),r') == d[{\rm els}(r)] :: C[{\rm els}(r')]}{
\begin{array}{l}
  \Gamma \rhd {\rm Sum}(A,B) \Longrightarrow C \; {\rm type} \\
  p: (\Gamma \rhd A \Longrightarrow d :: C [\text{Sum-sub-lf}(A,B)] \\
  q: (\Gamma \rhd B \Longrightarrow e :: C [\text{Sum-sub-rg}(A,B)] \\
  r:  (\Gamma \Longrightarrow a :: A) \\
  r':  (\Gamma \Longrightarrow  {\rm lf}(A,B,a,r):: {\rm Sum}(A,B)) \\
\end{array}}$$

$$\infer[(\text{Sum-c2})]{\Gamma \Longrightarrow \text{Sum-rec}(A,B,C,d,p,e,q,{\rm rg}(A,B,b,r),r') == e[{\rm els}(r)]:: C[{\rm els}(r')]}{
\begin{array}{l}
  \Gamma \rhd {\rm Sum}(A,B) \Longrightarrow C \; {\rm type} \\
  p: (\Gamma \rhd A \Longrightarrow d :: C [\text{Sum-sub-lf}(A,B)] \\
  q: (\Gamma \rhd B \Longrightarrow e :: C [\text{Sum-sub-rg}(A,B)] \\
  r:  (\Gamma \Longrightarrow b :: B) \\
  r':  (\Gamma \Longrightarrow  {\rm rg}(A,B,b,r):: {\rm Sum}(A,B)) \\
\end{array}}$$

$$\infer[(\text{Sum-cong})]{\Gamma \Longrightarrow{\rm Sum}(A,B) == {\rm Sum}(A',B')}{
\Gamma \Longrightarrow A == A' &  \Gamma \Longrightarrow B == B' }$$

$$\infer[(\text{lf-cong})]{\Gamma \Longrightarrow {\rm lf}(A,B,a,p) == {\rm lf}(A,B,a',p') :: {\rm Sum}(A,B) }{\Gamma \Longrightarrow B \; {\rm type} & p : (\Gamma \Longrightarrow a :: A) & p' : (\Gamma \Longrightarrow a' :: A) &\Gamma \Longrightarrow a == a' :: A }$$

$$\infer[(\text{rg-cong})]{\Gamma \Longrightarrow {\rm rg}(A,B,b,p) == {\rm rg}(A,B,b',p') :: {\rm Sum}(A,B) }{\Gamma \Longrightarrow A \; {\rm type}  & p : (\Gamma \Longrightarrow b :: B) & p' : (\Gamma \Longrightarrow b' :: B) &\Gamma \Longrightarrow b == b' :: B }$$

{\small
$$\infer[(\text{Sum-rec-cong})]{\Gamma \Longrightarrow \text{Sum-rec}(A,B,C,d,p,e,q,c,r) == \text{Sum-rec}(A',B',C',d',p',e',q',c',r') :: C[{\rm els}(r)]}{
\begin{array}{l}
  \Gamma \rhd {\rm Sum}(A,B) \Longrightarrow C \; {\rm type} \\
  \Gamma \rhd {\rm Sum}(A',B') \Longrightarrow C' \; {\rm type} \\
  p: (\Gamma \rhd {\rm Sum}(A,B) \Longrightarrow d :: C [\text{Sum-sub-lf}(A,B)] \\
  p': (\Gamma \rhd {\rm Sum}(A',B') \Longrightarrow d' :: C' [\text{Sum-sub-lf}(A',B')] \\
  q: (\Gamma \rhd {\rm Sum}(A,B) \Longrightarrow e :: C [\text{Sum-sub-rg}(A,B)] \\
  q': (\Gamma \rhd {\rm Sum}(A',B') \Longrightarrow e' :: C' [\text{Sum-sub-rg}(A',B')] \\
  r:  (\Gamma \Longrightarrow c:: {\rm Sum}(A,B)) \\
  r':  (\Gamma \Longrightarrow c':: {\rm Sum}(A',B')) \\
  Aq: (\Gamma \Longrightarrow A == A') \\
  Bq : (\Gamma \Longrightarrow B == B') \\
  \Gamma \rhd {\rm Sum}(A,B) \Longrightarrow \\ 
  \quad C == 
   C' [ \phi(\text{ext-eq'}({\rm Sum}(A, B), {\rm Sum}(A', B'), \text{Sum-cong}(A, A', B, B', Aq, Bq)) ]\\
   \Gamma \rhd A \Longrightarrow d == d' [ \phi(\text{ext-eq'}(A, A', Aq)) ] :: C [ \text{Sum-sub-lf}(A, B) ] \\
    \Gamma \rhd B \Longrightarrow e == e' [\phi(\text{ext-eq'}(B, B', Bq)) ] :: C [ \text{Sum-sub-rg}(A, B) ] \\
     \Gamma \Longrightarrow c == c' ::{\rm Sum}(A, B) 
\end{array}}$$
}

$$\infer[(\text{Sum-sub})]{\Delta \Longrightarrow{\rm Sum}(A,B)[h] == {\rm Sum}(A[h],B[h])}{
h: \Delta \longrightarrow \Gamma & \Gamma \Longrightarrow A \; {\rm type} &  \Gamma \Longrightarrow B \; {\rm type} }$$

$$\infer[(\text{lf-sub})]{\Gamma \Longrightarrow {\rm lf}(A,B,a,p) [h] ==  {\rm lf}(A[h],B[h],a[h],p') [h]  :: {\rm Sum}(A,B) [h]}{h: \Delta \longrightarrow \Gamma & \Gamma \Longrightarrow A \; {\rm type} & \Gamma \Longrightarrow B \; {\rm type} & p : (\Gamma \Longrightarrow a :: A)  & p' : (\Gamma \Longrightarrow a [h] :: A [h]) }$$

$$\infer[(\text{rg-sub})]{\Gamma \Longrightarrow {\rm rg}(A,B,b,p) [h] ==  {\rm rg}(A[h],B[h],b[h],p') [h]  :: {\rm Sum}(A,B) [h]}{h: \Delta \longrightarrow \Gamma & \Gamma \Longrightarrow A \; {\rm type} & \Gamma \Longrightarrow B \; {\rm type} & p : (\Gamma \Longrightarrow b :: B)  & p' : (\Gamma \Longrightarrow b [h] :: B[h]) }$$

{\small

$$\infer[(\text{Sum-rec-sub})]{
\begin{array}{l}
\Delta \Longrightarrow \text{Sum-rec}(A,B,C,d,p,e,q,c,r)[h]  \\
\qquad  == \text{Sum-rec}(A[h], B[h], C[\uparrow({\rm Sum}(A,B),h)], \\
\quad \qquad \qquad d[\uparrow(A,h)],p',e[\uparrow(B,h)],q', c[h],r') :: C[{\rm els}(r)]
\end{array}}{
\begin{array}{l}
 h : \Delta \longrightarrow \Gamma \\
  \Gamma \rhd {\rm Sum}(A,B) \Longrightarrow C \; {\rm type} \\
  p: (\Gamma \rhd {\rm Sum}(A,B) \Longrightarrow d :: C [\text{Sum-sub-lf}(A,B)]) \\
  q: (\Gamma \rhd {\rm Sum}(A,B) \Longrightarrow e :: C [\text{Sum-sub-rg}(A,B)])\\
  p': (\Delta \rhd (A[h]) \Longrightarrow  d[\uparrow(A,h)] :: C[\uparrow({\rm Sum}(A,B),h)][\text{Sum-sub-lf}(A[h],B[h])] )\\
  q': (\Delta \rhd (B[h]) \Longrightarrow  e[\uparrow(B,h)] ::  C[\uparrow({\rm Sum}(A,B),h)][\text{Sum-sub-rg}(A[h],B[h])])\\
  r:  (\Gamma \Longrightarrow c:: {\rm Sum}(A,B)) \\
  r':  (\Gamma \Longrightarrow c [h] :: {\rm Sum}(A[h],B[h])) \\
\end{array}}$$
}

\subsubsection{Universe rules}

For each $k \in {\mathbb N}$

$$\infer[(\text{U-}k)]{\Gamma \Longrightarrow {\rm U}_k \; {\rm type}}{\Gamma \; {\rm context}}
\qquad \infer[(TO DO)]{\Gamma \Longrightarrow A \; {\rm type}}{\Gamma \Longrightarrow A ::{\rm U}_k}$$

$$\infer[(\text{U-nat-})]{\Gamma \Longrightarrow {\rm Nat} :: {\rm U}_k }{\Gamma \; {\rm context}}
\qquad \infer[(\text{U-N0-})]{\Gamma \Longrightarrow {\rm N}_0 :: {\rm U}_k }{\Gamma \; {\rm context}}
$$

$$\infer[(\text{U-pi-})]{\Gamma \Longrightarrow \Pi_{\rm f}(A,B) :: {\rm U}_k }{\Gamma \Longrightarrow A ::  {\rm U}_k & \Gamma \rhd A \Longrightarrow B ::  {\rm U}_k }
\qquad
\infer[(\text{U-sigma-})]{\Gamma \Longrightarrow \Sigma_{\rm f}(A,B) :: {\rm U}_k }{\Gamma \Longrightarrow A ::  {\rm U}_k & \Gamma \rhd A \Longrightarrow B ::  {\rm U}_k }$$

$$\infer[(\text{U-Sum-})]{\Gamma \Longrightarrow {\rm Sum}(A,B) :: {\rm U}_k}{ \Gamma \Longrightarrow A ::  {\rm U}_k &
\Gamma \Longrightarrow B ::  {\rm U}_k }$$

$$\infer[(\text{U-ID-})]{\Gamma \Longrightarrow {\rm ID}(A,a,p,b,q) :: {\rm U}_k}{ \Gamma \Longrightarrow A ::  {\rm U}_k & p : (\Gamma \Longrightarrow a :: A) & q : (\Gamma \Longrightarrow b :: A)}$$

$$\infer[(\text{Cu-1a-})]{\Gamma \Longrightarrow {\rm U}_k :: {\rm U}_{s(k)}}{ \Gamma \; {\rm context}}
\qquad
\infer[(\text{Cu-1b-})]{\Gamma \Longrightarrow A :: {\rm U}_{s(k)}}{\Gamma \Longrightarrow A :: {\rm U}_k}$$

$$\infer[(\text{U-sub-})]{\Delta \Longrightarrow {\rm U}_k  [ h  ] ==  {\rm U}_k  }{h: \Delta \longrightarrow \Gamma}$$


$$\infer[(\text{U-eq-refl1})]{\Gamma \Longrightarrow  A ==  B  ::   {\rm U}_k }{\Gamma \Longrightarrow A ::  {\rm U}_k & \Gamma \Longrightarrow B ::  {\rm U}_k & \Gamma \Longrightarrow  A ==  B}$$

$$\infer[(\text{U-eq-refl2})]{\Gamma \Longrightarrow A == B \; {\rm type}}{\Gamma \Longrightarrow A  == B::{\rm U}_k}$$

\subsubsection{Bracket type rules} \label{bracket}

$$\infer{\Gamma \Longrightarrow{\rm Br}(A) \; {\rm type}}{
\Gamma \Longrightarrow A  \; {\rm type} } \qquad
\infer[(\text{Br-intro})]{\Gamma \Longrightarrow {\rm br}(a) :: {\rm Br}(A)}{
\Gamma \Longrightarrow a :: A  }$$

$$\infer[(\text{Br-e})]{\Gamma \Longrightarrow {\rm wh}(A,B,k,b,q,r,p) :: B}{
\begin{array}{l}
 \Gamma \Longrightarrow A  \; {\rm type} \\
  \Gamma \Longrightarrow B  \; {\rm type} \\
 q : (\Gamma \Longrightarrow k :: {\rm Br}(A)) \\
 r : (\Gamma \rhd A \Longrightarrow  b :: B [\downarrow(A)]) \\
 p : (\Gamma \rhd A  \rhd (A [\downarrow(A)])\Longrightarrow b[\text{pr-x}(A)] == b[\text{pr-y}(A)] ::  
 B [\downarrow(A)] [\downarrow (A [\downarrow(A)])])
 \end{array}}$$

$$\infer[(\text{Br-beta})]{\Gamma \Longrightarrow {\rm wh}(A,B,{\rm br}(a),b,q,r,p) == b[{\rm els}(t)]:: B}{
\begin{array}{l}
  \Gamma \Longrightarrow B  \; {\rm type} \\ 
 t : (\Gamma \Longrightarrow a :: A) \\
 q : (\Gamma \Longrightarrow k :: {\rm Br}(A)) \\
 r : (\Gamma \rhd A \Longrightarrow  b :: B [\downarrow(A)]) \\
 p : (\Gamma \rhd A  \rhd (A [\downarrow(A)])\Longrightarrow b[\text{pr-x}(A)] == b[\text{pr-y}(A)] ::  
 B [\downarrow(A)] [\downarrow (A [\downarrow(A)])])
 \end{array}}$$

$$\infer[(\text{Br-eta})]{\Gamma \Longrightarrow {\rm wh}(A,B,k,b [\text{br-sb}(A)],q,t,\text{br-sb-lm}(A, B, b,r)) == b[{\rm els}(q)]:: B}{
\begin{array}{l}
  \Gamma \Longrightarrow B  \; {\rm type} \\ 
 q : (\Gamma \Longrightarrow k :: {\rm Br}(A)) \\
 r : (\Gamma \rhd {\rm Br}(A) \Longrightarrow  b :: B [\downarrow({\rm Br}(A))]) \\
t : (\Gamma \rhd A \Longrightarrow b [\text{br-sb}(A)] :: B [\downarrow(A)]) 
 \end{array}}$$

$$\infer[(\text{Br-eqty})]{\Gamma \Longrightarrow a == b :: {\rm Br}(A)}{
\Gamma \Longrightarrow a :: {\rm Br}(A)  & \Gamma \Longrightarrow b :: {\rm Br}(A)}$$

$$\infer[(\text{Br-cong})]{\Gamma \Longrightarrow{\rm Br}(A) == {\rm Br}(A')}{
\Gamma \Longrightarrow A == A'}$$

$$\infer[(\text{Br-e-cong})]{\Gamma \Longrightarrow {\rm wh}(A,B,k,b,q,r,p) == {\rm wh}(A',B',k',b',q',r',p') :: B}
{\begin{array}{l}
 q : (\Gamma \Longrightarrow k :: {\rm Br}(A)) \\
 r : (\Gamma \rhd A \Longrightarrow  b :: B [\downarrow(A)]) \\
 p : (\Gamma \rhd A  \rhd (A [\downarrow(A)])\Longrightarrow b[\text{pr-x}(A)] == b[\text{pr-y}(A)] ::  
 B [\downarrow(A)] [\downarrow (A [\downarrow(A)])]) \\
  q' : (\Gamma \Longrightarrow k' :: {\rm Br}(A')) \\
 r' : (\Gamma \rhd A' \Longrightarrow  b' :: B' [\downarrow(A')]) \\
 p' : (\Gamma \rhd A'  \rhd (A' [\downarrow(A')])\Longrightarrow b'[\text{pr-x}(A')] == b'[\text{pr-y}(A')] ::  
 B' [\downarrow(A')] [\downarrow (A' [\downarrow(A')])]) \\
  Aq : (\Gamma \Longrightarrow  A == A') \\
  \Gamma \Longrightarrow  B == B' \\
  \Gamma \Longrightarrow k ==  k' :: {\rm Br}(A) \\
   \Gamma \rhd A \Longrightarrow  b == b' [ \phi (\text{ext-eq'}(A, A', Aq)) ] ::  B [\downarrow(A)] \\
 \end{array}}$$

$$\infer[(\text{Br-sub})]{\Delta \Longrightarrow {\rm Br}(A) [ f ]== {\rm Br}(A [f])}{
f: \Delta \longrightarrow \Gamma & \Gamma \Longrightarrow  A \; {\rm type}}$$

$$\infer[(\text{br-sub})]{\Delta \Longrightarrow {\rm br}(a) [ f ]== {\rm br}(a [f]) :: {\rm Br}(A) [ f ]}{
f: \Delta \longrightarrow \Gamma & \Gamma \Longrightarrow a :: A}$$

$$\infer[(\text{Br-e-sub })]{\Delta \Longrightarrow{\rm wh}(A, B, k, b, q, r, p) [ h ] == 
{\rm wh}(A [ h ], B [ h ], k [ h ], b [ \uparrow(A, h) ], q', r', p') :: B [h ]} 
{\begin{array}{l}
h : \Delta \longrightarrow \Gamma \\
 \Gamma \Longrightarrow A  \; {\rm type} \\
  \Gamma \Longrightarrow B  \; {\rm type} \\
 q : (\Gamma \Longrightarrow k :: {\rm Br}(A)) \\
 r : (\Gamma \rhd A \Longrightarrow  b :: B [\downarrow(A)]) \\
 p : (\Gamma \rhd A  \rhd (A [\downarrow(A)])\Longrightarrow b[\text{pr-x}(A)] == b[\text{pr-y}(A)] ::  
   B [\downarrow(A)] [\downarrow (A [\downarrow(A)])]) \\
 q' : (\Delta \Longrightarrow k[h] :: {\rm Br}(A[h])) \\
 r' : (\Delta \rhd A[h] \Longrightarrow  b[\uparrow(A,h)] :: B[h] [\downarrow(A[h])]) \\
 p' : (\Delta \rhd A[h] \rhd (A[h] [\downarrow(A[h])])\Longrightarrow b[\uparrow(A,h)][\text{pr-x}(A[h])] 
 == b[\uparrow(A,h)][\text{pr-y}(A[h])] \\   
  \qquad \qquad  :: B[h] [\downarrow(A[h])] [\downarrow (A[h] [\downarrow(A[h])])]) \\
 \end{array}}$$
 
 The universes are closed under bracket types

$$\infer[(\text{U-br-})]{\Gamma \Longrightarrow{\rm Br}(A) :: {\rm U}_k}{\Gamma \Longrightarrow A :: {\rm U}_k} $$

\subsection{Hidden arguments}

In the actual verification of the rules in Agda (Section \ref{overview}) there are implicit arguments that we have hidden in the above listing of rules. For instance 
$$\infer[(\Pi\text{-e})]{\Gamma \Longrightarrow \app(A,B,c,p,a,q) :: B[{\rm els}(q)]}{\Gamma \Longrightarrow A \; {\rm type} & \Gamma \rhd A  \Longrightarrow B \; {\rm type} & 
p : (\Gamma \Longrightarrow c :: \Pi_{\rm f}(A,B)) & q : (\Gamma \Longrightarrow a :: A)
}$$
looks as follows in Agda code:

\begin{verbatim}
Π-e :  {Γ : ctx} 
    -> (A : ty Γ)   
    -> (B : ty (Γ ▷ A))
    -> (c : raw Γ)
    -> (p : Γ ==> c :: Π-f {Γ} A B)
    -> (a : raw Γ)
    -> (q : Γ ==> a :: A)
--  -----------------------------------------
    ->  Γ ==> app A B c p a q :: B [[ els q ]]
\end{verbatim}

\section{Formalization in Agda} \label{overview}

This paper describes the formalization available at
 
 \medskip
 {\tt http://staff.math.su.se/palmgren/MLTT-and-setoids-2019-09-01.zip}.

 {\flushleft Needless to} say many improvements are possible, and some clean-ups of the code are probably necessary.
  Later improved or extended versions may be found at
 
 \medskip
 {\tt http://staff.math.su.se/palmgren/MLTT-and-setoids-latest.zip}. 
 
  \medskip
   The following files are part of the formalization we describe. Loading {\tt V-model-all-rules.agda} in Agda verifies all relevant files.  Agda version 2.5.2 has been used.

\medskip
Basic definitions and results concerning setoids

\medskip
{\tt

basic-types.agda

basic-setoids.agda      

dependent-setoids.agda   

subsetoids.agda   
}       

\medskip
Basic constructions and results concerning Aczel's iterative sets

\medskip
{\tt iterative-sets.agda      

iterative-sets-pt2.agda  

iterative-sets-pt3.agda  

iterative-sets-pt4.agda

iterative-sets-pt5.agda 

iterative-sets-pt6.agda

iterative-sets-pt8.agda
}

\medskip
The setoid model of extensional Martin-Löf type theory

\medskip
{\tt 

V-model-pt0.agda

V-model-pt1.agda         

V-model-pt2.agda        

V-model-pt3.agda         

V-model-pt4.agda         

V-model-pt5.agda         

V-model-pt6.agda         

V-model-pt7.agda         

V-model-pt8.agda         

V-model-pt9.agda         

V-model-pt10.agda

V-model-pt13.agda

V-model-pt11.agda

V-model-pt15.agda

V-model-all-rules.agda

}

\section{Comparing the Logical Framework and Agda} \label{LFAgda}

The Logical Framework (LF) is a dependently type lambda calculus which was designed to present dependent type theories in a compact form; see \cite{nordstrom+peterson+smith}. There is one basic type former for dependent products (or dependent function space)
$$\infer{\Gamma \Longrightarrow (x : \alpha)  \beta \; {\rm type}}{\Gamma, x : \alpha \Longrightarrow \beta \; {\rm type}}$$
It has an introduction rule which gives the only abstraction construction for terms, together with an elimination rules which is application.
$$\infer{\Gamma \Longrightarrow (x)b: (x : \alpha)  \beta}{\Gamma, x : \alpha \Longrightarrow b : \beta}
\qquad \infer{\Gamma \Longrightarrow c(a): \beta[a/x]}{\Gamma \Longrightarrow c : (x : \alpha)  \beta&  \Gamma \Longrightarrow a : \alpha }$$
There are corresponding $\beta$- and $\eta$-rules. Usual syntactic conventions are used to reduce the number of parentheses $c(a_1) \cdots(a_n)$ is abbreviated as $c(a_1,\ldots,a_n)$. If a type $\beta$ does not depend on $x$, $(x : \alpha)  \beta$ is abbreviated as $\alpha \rightarrow \beta$.
LF has one basic dependent type which is a type universe ${\rm Set}$ with a decoding function ${\rm El}(\cdot)$.
$$\infer{\Gamma \Longrightarrow {\rm Set} \; {\rm type}}{} \qquad 
\infer{\Gamma \Longrightarrow {\rm El}(a) \; {\rm type}}{\Gamma \Longrightarrow a : {\rm Set}}$$
A type theory T can now be axiomatized by introducing a number of new constants $c_1,\ldots,c_m$ with types in contexts
$$\infer{\Gamma_1\Longrightarrow c_1 : \alpha_1}{} \qquad \cdots \qquad  
\infer{\Gamma_m\Longrightarrow c_m : \alpha_m}{}$$
and furthermore equations in contexts
$$\infer{\Delta_1\Longrightarrow s_1 = t_1 : \beta_1}{} \qquad \cdots \qquad  
\infer{\Delta_n\Longrightarrow s_n = t_n: \beta_m}{}$$
In standard type theories the constants are type formers, introduction- and elimination-constants, and the equations express the computation rules. We refer to \cite{nordstrom+peterson+smith, ranta} for elaborations of Martin-Löf type theory in this form.

\medskip
The interactive proof system Agda has similar basic constructions (cf.\ right column)
$$\begin{array}{ll}
(x)b      & \backslash x \rightarrow b \\
(x : A)B & (x  :  A) \rightarrow B \\
a(b) & a\; b  \\
a(b,c) & a\; b \; c
\end{array}
$$
Agda has an infinite cumulative hierarchy of type universes ${\rm Set} = {\rm Set}0,   {\rm Set}1,  {\rm Set}2, \ldots $ with the rules
$$\infer{\Gamma \Longrightarrow {\rm Set}N \; {\rm type}}{} \qquad
\infer{\Gamma \Longrightarrow A: {\rm Set}(N+1)}{\Gamma \Longrightarrow A: {\rm Set}N} \qquad 
\infer{\Gamma \Longrightarrow A \; {\rm type}}{\Gamma \Longrightarrow A : {\rm Set}N}$$
Note that there is no explicit decoding function ${\rm El}$. In fact, every type in the system belongs to some ${\rm Set}N$ for some index $N$. Each universe ${\rm Set}N$ is closed under inductive-recursive definitions, which includes record types (generalized $\Sigma$ types) and recursive data types, as well as inductive families.

\section{Conclusion}

In this paper we have given a model of Martin-Löf extensional type theory with an infinite hiearchy 
of universes (cf.\ \cite{martinlof1982}) inside Martin-Löf intensional type theory (cf.\ \cite{nordstrom+peterson+smith}). The model is completely formalized in Agda using setoid constructions.
One may consider this as the first setoid model of the full extensional type theory in the intensional type theory.

\dontshow{
\section{Symmetric rules}

An alternative presentation of type theory is obtained by making the rules symmetric.
Basic judgement forms are then:

\begin{equation} \label{newinterpS}
\begin{array}{llll}
&\Gamma  \; {\rm context}  &\qquad &   \Gamma : {\mathbb V} \\
&(\Gamma \vdash A) \; {\rm type}  &\qquad &   A : [\kappa(\Gamma) \rightarrow {\mathbb V}]\\
&(\Gamma \vdash a)  \; {\rm raw} && a : [\kappa(\Gamma) \rightarrow {\mathbb V}] \\
\end{array}
\end{equation}
Then for these introduce equality judgements
\begin{equation} \label{newinterpS2}
\begin{array}{llll}
&\Gamma   == \Delta &\qquad &    \\
&(\Gamma \vdash A) == (\Delta \vdash B)  &\qquad &  \\
&(\Gamma \vdash a) == (\Delta \vdash b) &&  \\
\end{array}
\end{equation}
and a typing judgement:
\begin{equation} \label{newinterpS3}
\begin{array}{llll}
&(\Gamma \vdash a) ::  (\Delta \vdash A) &&  \\
\end{array}
\end{equation}
}

\bibliography{typesbib}{}
\bibliographystyle{plain}
\end{document}